\documentclass[11pt]{article}
\usepackage{mathrsfs}
\usepackage{amssymb}
\usepackage{amsmath}
\usepackage{amsbsy}
\usepackage{epsfig}
\usepackage{enumerate}
\usepackage{bm}
 \usepackage{color}

\usepackage[colorlinks,linkcolor=blue]{hyperref}

\newcommand{\MR}[1]{\href{http://www.ams.org/mathscinet-getitem?mr=#1}{\textcolor{blue}{MR-#1}}}%

\def\proclaim#1{\par \smallskip\noindent {\bf #1}\bgroup\it\ }
\def\endproclaim{\egroup\par\smallskip}

\newtheorem{lemma}{Lemma}[section]
\newtheorem{theorem}{Theorem}[section]
\newtheorem{proposition}{Proposition}[section]

\newtheorem{corollary}{Corollary}[section]

\newtheorem{remark}{Remark}[section]

\newbox\TempBox \newbox\TempBoxA

\def\pr{\textsf{P}} 
\def\ep{\textsf{E}} 

\def\Cov{\textsf{Cov}} 
\def\Var{\textsf{Var}} 

\def\underwiggle 1{
\ifmmode\setbox\TempBox=\hbox{$ 1$}\else\setbox\TempBox=\hbox{
1}\fi \setbox\TempBoxA=\hbox to \wd\TempBox{\hss\char'176\hss}
\rlap{\copy\TempBox}\smash{\lower9pt\hbox{\copy\TempBoxA}} }

\begin{document}

\title{\bf A stochastic algorithm approach for the elephant random walk with applications}

\author{ Li-Xin   ZHANG\footnote{Address: School of Statistics and Mathematics, Zhejiang Gongshang University, 310018, P. R. China; Center for Data Science,
Zhejiang University, 866 Yuhangtang Rd,
Hangzhou, 310058, P.R. China.  Email:stazlx@mail.zjgsu.edu.cn; stazlx\@ zju.edu.cn}\\
\\
 Zhejiang Gongshang University and
Zhejiang University, P.R. China}

\date{}
\maketitle
\begin{abstract}
  The randomized play-the-winner  rule (RPW) is a response-adaptive design proposed by Wei and Durham (1978) for sequentially randomizing patients to treatments in a two-treatment clinical trial so that more patients are assigned to the better treatment as the clinical trial goes on.  The elephant random walk (ERW)  proposed by Schutz and Trimper (2004) is a non-Markovian discrete-time random walk on $\mathbb Z$ which has a link to a famous saying that elephants can always remember where they have been. The asymptotic behaviors of RPW rule and  ERW have been studied in litterateurs independently, and their asymptotic behaviors are very similar. In this paper, we link RPW rule and ERW with the recursive stochastic algorithm. With the help of a recursive stochastic algorithm, we obtain the Gaussian approximation of the ERW as well as multi-dimensional varying-memory ERW with random step sizes.  By the Gaussian approximation, the central limit theorem, precise law of the iterated logarithm, and almost sure central limit theorem are obtained for  the multi-dimensional ERW, the multi-dimensional ERW with random step sizes, and their centers of   mass for all the diffusive, critical, and superdiffusive regimes. 
 Based on the Gaussian approximation and the small ball probabilities for a new kind of Gaussian process, the precise  Chung type laws of the iterated logarithm of the multi-dimensional ERW with random step sizes and its mass of center are also obtained for both the diffusive regime and superdiffusive regime.  

{\bf Keywords:}   elephant random walk, randomized paly-winner rule, recursive stochastic algorithm, Gaussian approximation, center of mass, small ball probability

{\bf AMS 2020 subject classifications:} Primary 60F15; Secondary 60F05

\end{abstract}

\section{Elephant random walk}
\setcounter{equation}{0}

The elephant random walk (ERW)  was introduced  by
Schütz and Trimper \cite{SchutaTrimper2004} in order to study the memory effects of a non-Markovian random
walk.  The model can be defined as follows. The elephant starts at time $n = 0$, with position $S_0=0$. At
time $n=1$, the elephant moves to $1$ with probability $p$ and to $-1$ with probability $1-p$,
where $p \in (0, 1]$. So the position of the elephant at time $n = 1$ is given by $S_1 =\sigma_1$,
with $\sigma_1$ a Rademacher $R(p)$ random variable. At time $n + 1$, $n \ge 1$, the step $\sigma_{n+1}$ is
determined stochastically by the following rule. Let $\beta_n$ be an integer which is chosen from
the set $\{1, 2,\ldots, n\}$ uniformly at random. If $\sigma_{\beta_n} = 1$, then $\sigma_{n+1}=1$ and $-1$ with probability $p$ and $1-p$, respectively. 
If $\sigma_{\beta_n} = -1$, then   $\sigma_{n+1}=1$ and $-1$ with probability $1-p$ and $p$, respectively. 
Equivalently, $\sigma_{n+1}$ is determined stochastically by the following rule:
$$\sigma_{n+1} =
\begin{cases}
 \sigma_{\beta_n} & \text{ with probability } p, \\
  -\sigma_{\beta_n} & \text{ with probability } 1-p.
 \end{cases}
 $$
 Thus, for $n \ge 2$, the position of the elephant at time $n$ is
\begin{equation}\label{eq:ERW1.1} S_n = \sum_{i=1}^n \sigma_i, \text{
 where } \sigma_n = \alpha_n \sigma_{\beta_n} , 
 \end{equation}
with $\alpha_n$ has a Rademacher distribution $R(p)$, $p \in (0, 1]$, and $\beta_n$ is uniformly distributed over
the integers $\{1, 2,\ldots, n-1\}$.  $p$ is called the  memory parameter.

The description of the asymptotic behavior of the ERW has motivated many interesting
works.  There is a critical phenomenon for the asymptotic behavior. The ERW is
respectively called  diffusive, critical and superdiffusive according to $p \in (0, 3/4)$, $p = 3/4$ and $p \in(3/4, 1]$.
By showing the connection of the ERW model with P\'olya-type urns, Baur and Bertoin \cite{BaurBertoin2016}
obtained the functional form of the central limit theorem. Coletti, Gava and Sch\"utz \cite{CGS2017a,CGS2017b}
proved the central limit theorem (CLT) and a strong invariance principle for $p \in (0, 3/4]$ and a
law of large numbers for $p \in (0, 1)$. Moreover, they also showed that if $p \in (3/4, 1]$, then the
ERW converges to a non-degenerate random variable which is not normal. V\'azquez Guevara
\cite{Vazque2019} gave the almost sure CLT. Bercu \cite{Bercu2018} recovered the CLT via a martingale method. Bertoin
\cite{Bertoin2022} studied how memory impacts passages at the origin for the ERW. Bercu and Laulin \cite{BercuLaulin2021}
investigated the asymptotic behavior of the center of mass of the ERW. Hu and Feng \cite{HuFeng2022} and  Dedecker et al \cite{DFHM2023} gave
an enhanced strong invariance principle of Coletti, Gava and Sch\"utz \cite{CGS2017b}.  Hayashi, Oshiro
and Takei \cite{HOT2023} obtained the exact rate of moment convergence in the CLT for the ERW.

Among these results, the following central limit theorem can be found in Baur and Bertoin \cite{BaurBertoin2016} and Bercu \cite{Bercu2018}. 
\begin{theorem}\label{th:ERWCLT} (\cite{BaurBertoin2016,Bercu2018}) If $0\le p<3/4$, then
$$\frac{S_n}{\sqrt{n}}\overset{\mathscr{D}}\to N\left(0,\frac{1}{3-4p}\right); $$
If $p=3/4$, then
$$\frac{S_n}{\sqrt{n\log n}}\overset{\mathscr{D}}\to N\left(0,1\right); $$
If $p>3/4$, then
$$ \frac{S_n}{n^{2p-1}}\to \xi \; a.s., $$
where   $\xi$ is a non-degenrate random variable.
\end{theorem}

The following   strong invariance principle (Gaussian approximation) can be found in Coletti,  Gava   and    Sch\"utz \cite{CGS2017a} and  Hu and Feng \cite{HuFeng2022}. 

\begin{theorem}\label{th:ERWInva} (\cite{CGS2017a,HuFeng2022})  Let $\{S_n;n\ge 1\}$   be the ERW with $p\le 4/3$ and let $\{ W(t); t\ge 0\}$ be a standard Brownian motion. Then there exists a common probability space for $\{S_n;n\ge 1\}$ and $\{ W(t); t\ge 0\}$  so that the following holds:
\begin{itemize}
  \item[\rm (i)] If $p<3/4$, then
 $$  \sqrt{3-4p} S_n - n^{2p-1}W(n^{3-4p}) =o\left(\sqrt{n \log\log n}\right) \; a.s.,$$
 $$ \sup_{0\le t\le 1}\left|\sqrt{3-4p} S_{[nt]} - (nt)^{2p-1}W((nt)^{3-4p})\right|=o(\sqrt{n}) \; \text{in probability}.$$
  \item[\rm (ii)] If $p=3/4$, then
 $$ S_n -n^{1/2}W(\log n) =o\left(\sqrt{n\log n\log\log \log n}\right) \; a.s., $$
  $$\sup_{0\le t\le 1}\left|n^{-t/2}S_{[n^t]} - W(t\log n)\right|=o(\sqrt{\log n}) \; \text{in probability}.$$
\end{itemize}
\end{theorem} 
Remark 1 of  Dedecker et al \cite{DFHM2023} gave a more precise approximation rate of the strong invariance principle.  

\section{Randomized play-the-winner rule}
\setcounter{equation}{0}

The randomized-play-the-winner (RPW) rule was introduced by Wei and Durham \cite{WD78}  for   sequentially randomizing patients to treatments.  Assume there
are two treatments (say, $TA$ and $TB$), with dichotomous response (success and
failure). The patients arrive at the clinical trial sequentially. Consider an urn with two types of balls (white and black) which starts at $W_0\ge 0$ white balls and $B_0\ge 0$ black balls.   For the $i$th patient, we a draw ball from the urn.  When the urn is empty, we assume that a type of white ball is drawn with  a probability $p_0$.  If a white ball is drawn, the patient is assigned to
the treatment $TA$, and otherwise, the patient is assigned to the treatment $TB$. The
ball is then replaced in the urn and the patient's response is observed.  A success
on treatment $TA$ or a failure on treatment $TB$ generates a white ball to the urn;
a success on treatment $TB$ or a failure on treatment $TA$ generates a black ball to the
urn. After  $n$ generations,   the number of white balls in the urn is denoted by $W_n$, the number of patients assigned to treatment $TA$ is denoted by $N_{nA}$.   
 
 Let $p_A = P(success|TA)$, $p_B  = P(success|TB)$, $q_A = 1-p_A$ and $q_B = 1-p_B$.  Assume $q_A\cdot q_B\ne  0$.
 Write $\rho=p_A-q_B$, $\alpha_0=W_0+B_0$, $v=q_B/(q_A+q_B)$. Then
 \begin{align}\label{eq:RPW1}
  \ep[W_n-W_{n-1}|W_1,\ldots,W_{n-1}]= & p_A\frac{W_{n-1}}{\alpha_0+{n-1}}+q_B\left(1-\frac{W_{n-1}}{\alpha_0+{n-1}}\right) \nonumber\\
 =& \rho\frac{W_{n-1}}{\alpha_0+{n-1}}+q_B=  \rho\frac{W_{n-1}-(\alpha_0+{n-1})v}{\alpha_0+{n-1}}+v. 
 \end{align}
 It follows that $\ep[W_1]=W_0+\rho\frac{W_0}{W_0+B_0}+q_B$ when $W_0+B_0\ne 0$, $=p_0p_A+(1-p_0)q_B$ when $W_0+B_0=0$, and
 \begin{align}\label{eq:RPWmean}
 \ep[W_n]-(\alpha_0+n)v=\prod_{k=1}^{n-1}\left(1+\frac{\rho}{\alpha_0+k}\right)\big(\ep W_1-(\alpha_0+1)v\big)
 \sim c n^{\rho}.
 \end{align}

 Wei and Durham \cite{WD78} showed that
 $$ \frac{W_n}{W_0+B_0+n}\to \frac{q_B}{q_A+q_B} \;a.s. \text{ and } \; \frac{N_{nA}}{n}\to \frac{q_B}{q_A+q_B} \; a.s. $$
 So, the RPW rule assigns more penitents to the better treatment. When the central limit theorem is considered, the critical phenomenon also appears. The following central limit theorem can be found in 
 Rosenberger\cite{Rosenberger1992} and Hu and Rosenberger \cite[Example 4.1]{HR06}.
\begin{theorem}\label{th:RPWCLT} (\cite{Rosenberger1992,HR06})
If $p_A+p_B<3/2$, then
$$\frac{W_n-n v}{\sqrt{n}}\overset{\mathscr{D}}\to N\left(0,\frac{q_Aq_B}{(q_A+q_B)^2(2(q_A+q_B)-1)}\right); $$
If $p_A+p_B=3/2$, then
$$\frac{W_n-n v}{\sqrt{n\log n}}\overset{\mathscr{D}}\to N\left(0,\frac{q_Aq_B}{(q_A+q_B)^2}\right); $$
If $p_A+p_B>3/2$, then
$$ \frac{S_n-nv}{n^{p_A+p_B-1}}\to W \; a.s., $$
where   $W$ is a non-degenrate random variable.
\end{theorem}

Let $\rho=p_A+p_B-1=1-q_A-q_B$, $\sigma^2=\frac{q_Aq_B}{(q_A+q_B)^2}$.
The Gaussian approximation of the RPW rule is obtianed by Bai, Hu and Zhang \cite{BHZ02} (c.f. Theorems 2.5, 2.6 and Corollaries 4.1, 4.2 there). 
\begin{theorem}  \label{th:RPWInva}(\cite{BHZ02})
Let $\{ B(t); t\ge 0\}$ be a standard Brownian motion. Then there exists a common probability space for $\{W_n;n\ge 1\}$ and $\{ B(t); t\ge 0\}$  so that the following holds:
\begin{itemize}
  \item[\rm (i)] If $\rho<1/2$, then
\begin{equation} \label{eq:thRPWInva1} W_n-\ep W_n-\sigma G_n=o\left(n^{1/2-\delta}\right) \; a.s. \;\; \forall 0<\delta < (1/2-\rho)\wedge (1/4), 
\end{equation}
 where
 $$ G_t=t^{\rho}\int_0^t s^{-\rho} d B(s) $$
 and 
 $$\{G_t;t\ge 0\}\overset{\mathscr{D}}=\left\{\frac{1}{\sqrt{1-2\rho}}t^{\rho}B(t^{1-2\rho});t\ge 0\right\}$$
  \item[\rm (ii)] If $\rho=1/2$, then for some $0<\delta<1/2$,
  \begin{equation} \label{eq:thRPWInva1} W_n-\ep W_n-\sigma\widehat{G}_n=o\left(n^{1/2}(\log n)^{1/2-\delta}\right)  \; a.s., \end{equation}
 where
 $$ \widehat{G}_t=t^{\rho}\int_1^t s^{-\rho} d B(s) $$
 and 
 $$\{\widehat{G}_t;t> 0\}\overset{\mathscr{D}}=\{t^{1/2}B(\log t );t> 0\}.$$
\end{itemize}
\end{theorem} 

Compare Theorems \ref{th:ERWCLT}, \ref{th:ERWInva} with Theorems \ref{th:RPWCLT}, \ref{th:RPWInva},  one can found that the   elephant random walk  and the randomized play-the-Winner rule behave asymptotically similarly. Next, we show that the randomized play-the-Winner rule  can be regarded as  a biased elephant random walk.

\bigskip
{\em RPW rule as biased ERW}.
A biased ERW in which the elephant has  a favorite direction  is defined as follows. The elephant starts at time $n = 0$, with position $S_0=0$. At
time $n=1$, the elephant moves to $1$ with probability $p_A$ and to $-1$ with probability $1-p_A$. So the position of the elephant at time $n = 1$ is given by $S_1 =\sigma_1$,
with $\sigma_1$ a Rademacher $R(p_A)$ random variable. At time $n + 1$, $n \ge 1$, the step $\sigma_{n+1}$ is
determined stochastically by the following rule. Let $\beta_n$ be an integer which is chosen from
the set $\{1, 2,\ldots, n\}$ uniformly at random. If $\sigma_{\beta_n} = 1$, then   $\sigma_{n+1}=1$ and $-1$ with probability $p_A$ and $1-p_A$, respectively. If $\sigma_{\beta_n} = -1$, then   $\sigma_{n+1}=1$ and $-1$ with probability $1-p_B$ and $p_B$, respectively. 
When $p_A>p_B$, the elephant  prefers to walk to the right. 
For $n \ge 2$, the position of the elephant at time $n$ is
\begin{equation}\label{eq:ERW1.1} S_n = \sum_{i=1}^n \sigma_i. 
 \end{equation}
It is easily seen that for $n\ge 2$, 
$$  \ep[ \sigma_n| S_1,\ldots,S_{n-1}]=(p_A+p_B-1)\frac{S_{n-1}}{n-1}+p_A-p_B. $$

 Consider a play-the-winner  rule with $W_0=B_0=0$,  and to start, a white ball is drawn with probability $p_0=1$. Let $\widetilde{S}_n=2W_n-n$, $\widetilde{\sigma}_n=\widetilde{S}_n
-\widetilde{S}_{n-1}$. Then $\widetilde{S}_0=0$, $\widetilde{\sigma}_n=\pm 1$ and 
$$ \ep[\widetilde{\sigma}_n|\widetilde{S}_1,\ldots,\widetilde{S}_{n-1}]=(p_A+p_B-1)\frac{\widetilde{S}_{n-1}}{n-1}+p_A-p_B,\; n\ge 2, $$
by \eqref{eq:RPW1}. 
Thus
$$ \{S_n;n\ge 0\}\overset{\mathscr{D}}=\{ \widetilde{S}_n;n\ge 0\} \;\text{ when } W_0=B_0=0, p_0=1. $$

By Theorem \ref{th:RPWInva} and  \eqref{eq:RPWmean}, we have 
\begin{itemize}
  \item[\rm (i)] If $\rho=p_A+p_B-1<1/2$, then for any  $0<\delta < (1/2-\rho)\wedge (1/4)$, 
$$  S_n-\frac{p_A-p_B}{2-p_A-p_B}-\frac{2\sigma}{\sqrt{1-2\rho}} W(n^{1-2\rho})=o\left(n^{1/2-\delta}\right) \; a.s.   
$$
  \item[\rm (ii)] If $\rho=p_A+p_B-1=1/2$, then for some $0<\delta<1/2$,
 $$ S_n-\frac{p_A-p_B}{2-p_A-p_B}-2\sigma n^{1/2} W(\log n)=o\left(n^{1/2}(\log n)^{1/2-\delta}\right)  \; a.s.  $$
\end{itemize}
 where $\sigma^2=q_Aq_B/(q_A+q_B)^2$.
When $p_A=p_B=p$, $\{\widetilde{S}_n;n\ge 0\}$ is   the original ERW.


\section{Gaussian approximation by using  stochastic algorithm}
\setcounter{equation}{0}

Stochastic approximation  algorithms, which
have progressively gained sway thanks to the development of computer science and automatic control theory,  have been the subject of many studies. A   recursive stochastic approximation algorithm is defined on a filtered probability space $\big(\Omega,\mathscr{F}, (\mathscr{F}_n)_{n\ge 0}, \pr)$ and has the following structure:
\begin{equation}\label{eq:SAModel} \bm\theta_{n+1}=\bm\theta_n-\gamma_{n+1} \bm h(\bm\theta_n) + \gamma_{n+1}(\Delta \bm M_{n+1}+\bm r_{n+1}),
\end{equation}
where $\bm\theta_n$ is a row vector in $\mathbb R^d$, the regression function $\bm h: \mathbb R^d \to \mathbb R^d$ is a real vector-valued function, $\bm \theta_{0}$ is a finite random vector, $\bm M_0=\bm 0$, $\{\Delta \bm M_n,\mathscr{F}_n;n\ge 1\}$ is a sequence of  martingale differences and $\bm r_n$ is a residual term, $\{\gamma_n\}$ is a positive sequence of steps  that tends toward zero, such that $\sum_{n=1}^{\infty}\gamma_n$ diverges.  
The basic asymptotic properties of  the recursive stochastic algorithms   can be found in text books as  Kushner  and   Clark \cite{KC78},  Duflo \cite{D97}, Benveniste et al. \cite{BMP90} and Kushner and Yin \cite{KY03}   with different groups of conditions.

For the RPW rule, by \eqref{eq:RPW1},
$$ W_n-\ep W_n=(W_{n-1}-\ep W_{n-1})+\rho\frac{W_{n-1}-\ep W_{n-1}}{\alpha_0+n-1}+\Delta M_n, $$
where $\rho=p_A-q_B$, and $\Delta M_n=(W_n-W_{n-1})-\ep[W_n-W_{n-1}|W_1,\ldots,W_{n-1}]$ is a sequence of martingale differences. Thus, $\theta_n=(W_n-\ep W_n)/n$ satisfies the one-dimensional  recursive stochastic algorithm:
\begin{equation}\label{SA1}
\theta_n=\theta_{n-1}-  \frac{h(\theta_{n-1})}{n}+\frac{\Delta M_n+r_n}{n}
\end{equation}
 with the    regression function $h(x)=(1-\rho)x$ and residual $r_n=-\alpha_0W_{n-1}/\{(n-1)(\alpha_0+n-1)\}$.

For the ERW, it is easily seen that
\begin{equation}\label{eq:ERW1.2}
\ep\left[\sigma_n |S_1,\ldots,S_{n-1}\right]=\rho\frac{S_{n-1}}{n-1},
\end{equation}
where $\rho=2p-1$.  Thus, 
$$ S_n=S_{n-1}+\rho\frac{S_{n-1}}{n-1}+\Delta M_n, $$
where $\Delta M_n=:M_n-M_{n-1}=\sigma_n-\ep[\sigma_n|S_1,\ldots,S_{n-1}]$ is a sequence of martingale  differences.  
If let $\theta_n=S_n/n$, then $\theta_n$ satisfy the  
recursive stochastic algorithm \eqref{SA1} 
 with the regression function     $h(x)=(1-\rho)x$ the same as that of the RPW rule, and residual  $r_n\equiv 0$.
Since the martingale behaves asymptotically like a Brownian motion,    the ERW and RPW rule with the same recursive stochastic algorithm structure behave  asymptotically similarly. 
 
Laruelle and Pag\`es \cite{LP13} linked multi-color urn models with a multi-dimensional recursive stochastic algorithm.  Zhang \cite{Zhang2016} obtained the central limit theorem and Gaussian approximation of a general multi-dimensional recursive stochastic algorithm and studied the asymptotic properties of general randomized urn models by the recursive stochastic algorithm. The critical phenomenon of the asymptotic behaviors of the recursive stochastic algorithms is described by the smallest real part of the eigenvalues of the gradient matrix of the regression function $\bm h$. The central limit theorem and Gaussian approximation of the  ERW can be also derived from those for the recursive stochastic algorithm.    The main purpose of this paper is to use the theory for the recursive stochastic algorithm to obtain the Gaussian approximation of the multi-dimensional ERW  and the multi-dimensional ERW with random step sizes. By the Gaussian approximation,  many other asymptotic properties such as the central limit theorem,   law of the iterated logarithm, and almost sure central limit theorem of the ERWs with random step sizes and their centers of mass are obtained for all the diffusive, critical, and superdifficult regimes, by deriving the same properties of related Gaussian processes. Especially, the Chung type law of the iterated logarithm of the multi-dimensional ERW is established.   The method can also be used to study other types of elephant random walks, for example, 
the reinforced ERW model proposed by Laulin \cite{Laulin2022}, the ERW with stops \cite{Bercu2022}.

The multi-dimensional type ERW  was introduced  by Bercu and
 Laulin \cite{BercuLaulin2019}. It is a  discrete-time random walk on $\mathbb Z^d$ 
where $d\ge 1$. The elephant starts at the
origin at time zero, $\bm S_0 = \bm 0$. At time $n = 1$, the elephant moves in one of the $2d$ directions
with the same probability $1/(2d)$. Afterwards, at time $n + 1$, the elephant chooses uniformly
at random an integer $\beta_n$ among the previous times $1, \cdots, n$. Then, it moves exactly in the same
direction as that of time $\beta_n$ with probability $p$ or in one of the $2d-1$ remaining directions with
the same probability $(1-p)/(2d-1)$, where the parameter $p \in [0, 1]$ stands for the memory
parameter of the ERW.  Therefore, the position of the elephant at time $n+1$  is given
by
$$ \bm S_{n+1} = \bm S_n + \bm \sigma_{n+1}, $$
where the vector $\bm\sigma_{n+1}=(\sigma_{n+1,1},\ldots,\sigma_{n+1,d})$ is the $(n+1)$th increment of the random walk with one of its elements $\sigma_{n+1,k}$  is $\pm 1$, and others are $0$. The ERW shows three different
regimes, diffusive regime, critical regime and superdiffusive regime, depending on the location of its memory parameter $p$ with respect to the critical value, 
$$ p_d=\frac{2d+1}{4d}. $$  For the multi-dimensional ERW, the central limit theorem and weak invariance principle are obtained by Bercu and
Laulin \cite{BercuLaulin2019}, Bercu and
Laulin \cite{BercuLaulin2021} and  Bertenghi \cite{Bertenghi2022}. 

As an extension of the ERW, 
an ERW with random step sizes may be described  by
$$ \bm T_n=\sum_{k=1}^n \bm \sigma_kZ_k, $$
where $\{Z_n;n\ge 1\}$ is a sequence of independent random variables which are independent of $\{\bm S_n;n\ge 1\}$, $Z_i$ is the step size of the $i$th step. Assume that $\mu_Z=\ep[Z_1]$ and  $\sigma_Z^2=\Var(Z_1)$ are finite, and 
\begin{equation}\label{eq:momentcondition}\ep| Z_1|^{2+\epsilon}<\infty \text{ for some } \epsilon>0.
\end{equation}
  Then
$$\bm T_n=\sum_{k=1}^n \bm\sigma_k(Z_k-\mu_Z)+\mu_Z\bm S_n. $$
The one-dimensional ERW with random step sizes
was introduced by Fan and Shao \cite{FanShao2023}, and  Dedecker et al \cite{DFHM2023} obtained the law of the iterated logarithm, the central limit theorem and rates of convergence in the central limit theorem. 

In the above setting, the memory
parameter $p$ of the ERW is the same at each step. In general, we consider a varying-memory ERW with random steps sizes in which the memory parameter is $p_i$ at the $i$-th step. That is, at the $i$-th step, the elephant  moves exactly in the same
direction as one of the previous steps  with probability $p_i$ or in one of the $2d-1$ remaining directions with
the same probability $(1-p_i)/(2d-1)$. We assume that
\begin{equation}\label{eq:non-hom} \frac{\sum_{i=1}^n p_i}{n}-p=o(n^{-\epsilon}) \text{ for some } \epsilon>0. 
\end{equation}
We consider the joint asymptotic behavior of $\bm S_n$ and $\bm T_n$.  
Let 
$$\rho=\frac{2dp-1}{2d-1},\; \rho_i=\frac{2dp_i-1}{2d-1}, $$
 $\Sigma=\{\pm e_1,\ldots,\pm e_d\}$, 
where $e_k=(0,\ldots,0,1,0,\ldots)$ is the vector with the $k$th elements  $1$, and $0$ others. Then 
$$\pr(\bm\sigma_{n+1}=\bm\sigma_{\beta_n}|\beta_n)=p_{n+1}\bm\sigma_{\beta_n}, \;\; \pr(\bm\sigma_{n+1}=\bm\sigma|\beta_n)=\frac{1-p_{n+1}}{2d-1}, \bm\sigma\ne \bm\sigma_{\beta_n}. $$
Thus
\begin{align}\label{eq:condExp1}
 \ep[\bm\sigma_{n+1}|\beta_n]=&p_{n+1}\bm\sigma_{\beta_n}+\frac{1-p_{n+1}}{2d-1}\sum_{\bm\sigma\in \Sigma\setminus \bm\sigma_{\beta_n}}\bm\sigma  
 = p_{n+1}\bm\sigma_{\beta_n}+\frac{1-p_{n+1}}{2d-1}\left(\sum_{\bm\sigma\in \Sigma}\bm\sigma-\bm\sigma_{\beta_n}\right)\nonumber \\
 =&p_{n+1}\bm\sigma_{\beta_n}+\frac{1-p_{n+1}}{2d-1}\left(\bm 0-\bm\sigma_{\beta_n}\right)=\rho_{n+1} \bm\sigma_{\beta_n}.
 \end{align}
It follows that
$$\ep[\bm\sigma_{n+1}|\bm S_1,\ldots,\bm S_n]=\rho_{n+1} \frac{\bm S_n}{n}.$$
Thus
\begin{align}\label{eq:SAforS1}\bm S_{n+1}=&\bm S_n+\rho_{n+1} \frac{\bm S_n}{n}+\bm\sigma_{n+1}-\ep[\bm\sigma_{n+1}|\mathscr{F}_n] \\
\label{eq:SAforS2}
=&\bm S_n+\rho \frac{\bm S_n}{n}+\bm\sigma_{n+1}-\ep[\bm\sigma_{n+1}|\mathscr{F}_n]+(\rho_{n+1}-\rho) \frac{\bm S_n}{n},
\end{align}
where $\mathscr{F}_n=\sigma(\bm S_k,Z_k; k\ldots,n)$. Equivelently,
$$ \frac{\bm S_{n+1}}{n+1}=\frac{\bm S_n}{n}-\frac{(1-\rho)\bm S_n/n}{n+1}+\frac{\bm\sigma_{n+1}-\ep[\bm\sigma_{n+1}|\mathscr{F}_n]}{n+1}+\frac{(\rho_{n+1}-\rho)\bm S_n/n}{n+1}. $$
 Let 
\begin{equation} \label{eq:martingales.1}
\Delta\bm M_{n+1,1}=\bm \sigma_{n+1}-\ep[\bm\sigma_{n+1}|\mathscr{F}_n],  \;\; \Delta\bm M_{n+1,2}=\bm \sigma_{n+1}(Z_{n+1}-\mu_Z), \end{equation}
   Then $\bm \theta_n=\big(\bm S_n,\bm T_n)/n$  satisfies the recursive stochastic algorithm:
 $$\bm \theta_{n+1}=  \bm \theta_n-\frac{\bm h(\bm \theta_n)}{n+1}+\frac{(\Delta \bm M_{n+1,1}, \Delta \bm M_{n+1,2}+\mu_Z\Delta\bm M_{n+1,1})}{n+1}+\frac{\bm r_{n+1}}{n+1}, $$
 with
 \begin{align*}
 & \bm h(\bm x,\bm y)=\big((1-\rho)\bm x, \bm y-\mu_Z \bm x\big)=(\bm x,\bm y)\bm H, \\
 & D \bm h(\bm x,\bm y)=\bm H=\begin{pmatrix} 1-\rho  & -\mu_Z \\
   0 & 1\end{pmatrix}\otimes I_d.
 \end{align*}

By the theorems  on the recursive stochastic algorithm (c.f. Theorems 3.1  and 2.2 of Zhang \cite{Zhang2016}), we have the following Gaussian approximation. 
\begin{theorem}  \label{th:multi-RPWInva}
Assume \eqref{eq:momentcondition} and \eqref{eq:non-hom}. 
Let $\{ B_k(t), W_k(t); t\ge 0, k=1,\ldots,d\}$ be $2d$ independent  standard Brownian motions, $\bm W(t)=\big(W_1(t),\ldots,W_d(t)\big)$, $\bm B(t)=\big(B_1(t),\ldots,B_d(t)\big)$. Then there exists a common probability space and $\tau\in (0,1/2)$ for $\{\bm S_n, \bm T_n;n\ge 1\}$ and $\{ B_k(t), W_k(t); t\ge 0, k=1,\ldots,d\}$  so that the following holds:
\begin{itemize}
  \item[\rm (i)] If $\rho<1/2$ ($p<p_d$), then
\begin{align} \label{eq:multi-thRPWInva1}
&  \bm S_n=\frac{1}{\sqrt{d }} \bm G_n+o\big(n^{1/2-\tau}\big) \; a.s., \\
\label{eq:multi-thRPWInva2}
&\bm T_n=\frac{\sigma_Z}{\sqrt{d}}\bm W(n)+\frac{\mu_Z}{\sqrt{d }} \bm G_n+o\big(n^{1/2-\tau}\big) \; a.s.,  
\end{align}
 where  
  \begin{align} \label{eq:multi-thRPWInva3} \bm G_t=t^{\rho}\int_0^t s^{-\rho} d \bm B(s), \;\; \rho<1/2,
  \end{align}
  is a solution of the differential equation 
  \begin{equation}\label{eq:SDE}
  d\bm G_t=\rho \frac{\bm G_t}{t} dt +d \bm B(t) 
  \end{equation}
  with $\bm G_0=\bm 0$.
  \item[\rm (ii)] If $\rho=1/2$ ($p=p_d$), then  
 \begin{align} \label{eq:multi-thRPWInva4}
&  \bm S_n=\frac{1}{\sqrt{d }} \widehat{\bm G}_n+O(n^{1/2}) \; a.s., \\
\label{eq:multi-thRPWInva5}
&\bm T_n=\frac{\sigma_Z}{\sqrt{d}}\bm W(n)+\frac{\mu_Z}{\sqrt{d }} \widehat{\bm G}_n+O(n^{1/2}) \; a.s.,  
\end{align}
 where  
\begin{align} \label{eq:multi-thRPWInva6} \widehat{\bm G}_t=t^{\rho}\int_1^t s^{-\rho} d \bm B(s)
\end{align}
is a solution of \eqref{eq:SDE} with $\bm G_1=\bm 0$. 
  \item[\rm (iii)] If $\rho>1/2$ ($p>p_d$), then there is a  non-degenerate random vector $\bm \xi\in\mathbb R^d$ such that
   \begin{align} \label{eq:multi-thRPWInva7} \frac{\bm S_n}{n^{\rho}}\to \bm \xi \; a.s. \text{ and } \frac{\bm T_n}{n^{\rho}}\to \mu_Z\bm \xi\; a.s.
   \end{align}
   \end{itemize}
\end{theorem} 

For the   superdiffusive regime, we have a more precise result. 
\begin{theorem}  \label{th:multi-RPWInva2} Suppose $\rho>1/2$.
Instead of  \eqref{eq:momentcondition} and \eqref{eq:non-hom}, we assume $\ep Z_1^2<\infty$ and 
\begin{equation}\label{eq:non-hom2}
\sum_{i=1}^{\infty} \frac{p_i-p}{i} \text{ is convergent}. 
\end{equation}
 \begin{enumerate}
   \item[\rm (a)] Let $\digamma_n=n^{1/2-\rho}+|\sum_{i=n}^{\infty}(\rho_i-\rho)/i|$. Then \eqref{eq:multi-thRPWInva7} holds, and
 \begin{align} \label{eq:multi-thRPWInva8}  \left\|\frac{\bm S_n}{n^{\rho}}- \bm \xi \right\|_r=O(\digamma_n) \text{ for any } r>0 \text{ and }
   \left\|\frac{\bm T_n}{n^{\rho}}- \mu_Z\bm \xi \right\|_q=O(\digamma_n), 
  \end{align}
  if $\ep|Z_1|^q<\infty$, $q\ge 2$, where $\|X\|_r=(\ep|X|^r)^{1/r}$.
  \item[\rm (b)]  $\ep[\bm \xi]=0$ and
 \begin{align} \label{eq:multi-thRPWInva12}   \ep[\bm \xi^{\prime}\bm \xi]=\frac{C}{d} I_d, 
 \end{align}
  where $C>0$ is a constant which only depends on $\rho_i$s.
   \item[\rm (c)] If $1/2<\rho<1$ ($p_d<p<1$), then  for any nonzero  $\bm b\in\mathbb R$, $\bm\xi\bm b^{\prime}$ has   no point probability mass, namely, $\pr(\bm\xi\bm b^{\prime}=x) = 0$  for any $x$. 
   \item[\rm (d)] If   $1/2<\rho<1$ ($p_d<p<1$), \eqref{eq:momentcondition} holds, and \eqref{eq:non-hom2} is strengthened  to 
  \begin{equation}\label{eq:non-hom3}
\frac{1}{n}\sum_{i=1}^n p_i-p=O(n^{-\epsilon_0}), \;\; \epsilon_0>\rho-1/2, 
\end{equation}
then for some $\tau>0$,
\begin{align} \label{eq:multi-thRPWInva13} 
&\bm S_n-  n^{\rho}\bm \xi = \frac{1}{\sqrt{d}}\bm G_n +o(n^{1/2-\tau})\; a.s.,\\
&  \label{eq:multi-thRPWInva14} \bm T_n-n^{\rho} \mu_Z\bm \xi =\frac{\sigma_Z}{\sqrt{d}}\bm W(n)+\frac{\mu_Z}{\sqrt{d}}\bm G_n+o(n^{1/2-\tau})\; a.s.
\end{align}
where 
  \begin{align} \label{eq:multi-thRPWInva15} \bm G_t=-t^{\rho}\int_t^{\infty} s^{-\rho} d \bm B(s),\; \; \rho>1/2,
  \end{align}
  is a solution of \eqref{eq:SDE} with $\bm G_t=\bm 0$. 
 \end{enumerate}
 \end{theorem}

By Theorems \ref{th:multi-RPWInva}, \ref{th:multi-RPWInva2}, and the properties of Gaussian processes, we have the following central limit theorem and 
 law of the iterated logarithm for all the diffusive, critical and superdiffusive regimes. 
\begin{corollary} \label{cor:multi-ERW}Consider the ERW with random step sizes.  
Assume \eqref{eq:momentcondition} and \eqref{eq:non-hom}.
\begin{itemize}
  \item[\rm (i)] If $\rho<1/2$ ($p<p_d$), then
\begin{align} \label{eq:multi-thERWCLT1}
&\frac{(\bm S_n, \bm T_n)}{\sqrt{n}}\overset{\mathscr{D}}\to N(\bm 0,\bm \Sigma_{\rho}) \text{ with } \bm \Sigma_{\rho}= \begin{pmatrix} \frac{1}{1-2\rho} & \frac{\mu_Z}{1-2\rho} \\
\frac{\mu_Z}{1-2\rho} & \sigma_Z^2+\frac{\mu_Z^2}{1-2\rho}
\end{pmatrix}\otimes \frac{I_d}{d},
\end{align}
\begin{align} \label{eq:multi-thERWLIL1}
&\limsup_{n\to \infty} \frac{\|\bm T_n\|}{\sqrt{2n\log\log n}} =  \sqrt{\frac{\sigma_Z^2}{d}+\frac{\mu_Z^2}{d(1-2\rho)}}\; a.s. 
\end{align} 
  \item[\rm (ii)] If $\rho=1/2$ ($p=p_d$), then  
 \begin{align} \label{eq:multi-thERWCLT2}
 \frac{(\bm S_n, \bm T_n)}{\sqrt{n\log n}} \overset{\mathscr{D}}\to N \Big(\bm 0,  \widetilde{\bm \Sigma}\Big) \; \text{ with }   \; \widetilde{\bm \Sigma}=\begin{pmatrix} 1 &\mu_{Z} \\ \mu_Z & \mu_Z^2\end{pmatrix} \otimes  \frac{I_d}{d},
\end{align}
\begin{align} \label{eq:multi-thERWLIL2}
&\limsup_{n\to \infty} \frac{\|\bm T_n\|}{\sqrt{2n(\log n)(\log\log \log n)}}=  \frac{|\mu_Z| }{\sqrt{d}} \; a.s. 
\end{align}
  \end{itemize}
 \item[\rm (iii)] If $1/2<\rho<1$ ($p_d<p<1$) and \eqref{eq:non-hom3} holds, then
\begin{align} \label{eq:multi-thERWCLT3}
&\frac{(\bm S_n-n^{\rho}\bm \xi, \bm T_n-\mu_Zn^{\rho}\bm \xi)}{\sqrt{n}}\overset{\mathscr{D}}\to N(\bm 0,\bm \Sigma_{\rho}) \text{ with } \bm \Sigma_{\rho}= \begin{pmatrix} \frac{1}{2\rho-1} & \frac{\mu_Z}{2\rho-1} \\
\frac{\mu_Z}{2\rho-1} & \sigma_Z^2+\frac{\mu_Z^2}{2\rho-1}
\end{pmatrix}\otimes \frac{I_d}{d},
\end{align}
\begin{align} \label{eq:multi-thERWLIL3}
&\limsup_{n\to \infty} \frac{\|\bm T_n-\mu_Zn^{\rho}\bm \xi\|}{\sqrt{2n\log\log n}} =  \sqrt{\frac{\sigma_Z^2}{d}+\frac{\mu_Z^2}{d(2\rho-1)}}\; a.s. 
\end{align} 
\end{corollary}

\begin{remark} When $Z_i\equiv \mu_Z=1$, $\bm T_n=\bm S_n$. Thus, the central limit theorem and the law of the iterated logarithm are consistent with those established by   Bercu and Laulin \cite{BercuLaulin2019} for $\rho\le 1/2$. When $\rho>1/2$ and $p_i\equiv p$, 
Bercu and Laulin \cite{BercuLaulin2019} obtained the a.s. and $L_2$ convergence of $\bm S_n/n^{\rho}$, and showed that $\ep[\bm\xi]=0$,
$$ \ep[\bm\xi^{\prime}\bm \xi]=\frac{1}{d(2\rho-1)\Gamma(2\rho)} I_d. $$ 
\end{remark} 

\begin{remark}
For the one-dimensional case $(d=1)$, \eqref{eq:multi-thERWCLT1} and \eqref{eq:multi-thERWCLT2} are consistent with (3.1) of Dedecker et al \cite{DFHM2023}, and \eqref{eq:multi-thERWLIL2} is consistent with (2.2) of  Dedecker et al \cite{DFHM2023}. When $\rho<1/2$,
Dedecker et al \cite{DFHM2023} only gave an upper bound of the law of the iterated logarithm as 
$$ \limsup_{n\to \infty} \frac{| T_n|}{\sqrt{2n\log\log n}}\le \sigma_Z+\frac{|\mu_Z|}{\sqrt{1-2\rho}}\; a.s., \; \rho<1/2. $$
The upper bound is larger than the exact value on the right hand of \eqref{eq:multi-thERWLIL1}. 

In Dedecker et al \cite{DFHM2023}, only the second moment $\sigma_Z<\infty$ is assumed. By using the truncation method, in Corollary \ref{cor:multi-ERW},   the condition $\ep|Z_1|^{2+\epsilon}$ can be removed. In fact, if denote $Z_i^{(1)}=Z_iI\{|Z_i|\le c\}$, $Z_i^{(2)}=Z_iI\{|Z_i|> c\}$,
$$ \bm T_n^{(j)}=\sum_{k=1}^n \bm \sigma_kZ_k^{(j)}, \;\; j=1,2. $$
Then \eqref{eq:multi-thERWCLT1}-\eqref{eq:multi-thERWLIL3} hold for $\bm T_n^{(1)}$. Also, the same limit values for $\bm T_n^{(2)}$ will go to zero as $c\to \infty$. Take \eqref{eq:multi-thERWLIL1} as an example,
\begin{align*}  
&\limsup_{n\to \infty} \frac{\|\bm T_n^{(2)}\|}{\sqrt{2n\log\log n}}\\
\le & \limsup_{n\to \infty}\sum_{l=1}^d \frac{|\sum_{m=1}^n \sigma_{m,l} (Z_m^{(2)}-\ep[Z_m^{(2)}])|}{\sqrt{2n\log\log n}} 
+\big|\ep[Z_1^{(2)}]\big|\limsup_{n\to \infty} \frac{\|\bm S_n\|}{\sqrt{2n\log\log n}} \\
\le & d\sqrt{\Var(Z_1^{(2)})}+   \frac{\big|\ep[Z_1^{(2)}]\big|}{\sqrt{d(1-2\rho)}}\to 0 \; \text{ as } c\to \infty. 
\end{align*}
\end{remark}

The following is the Chung type law of the iterated logarithm of   $\bm T_n$. 
\begin{corollary} \label{cor:CLIL} \begin{itemize}
\item[(i)]Assume \eqref{eq:momentcondition} and \eqref{eq:non-hom}. If $\rho<1/2$ ($p<p_d$), then
\begin{equation}\label{eq:CLIL2}\liminf_{n\to \infty}\Big(\frac{\log\log n}{n}\Big)^{1/2} \max_{m\le n}\|\bm T_m\|=\frac{j_{\nu}}{\sqrt{2d}}(\ep[Z_1^2])^{1/2} \; a.s., 
\end{equation} 
where $j_{\nu}$ denotes the smallest positive zero of the Bessel function $J_{\nu}(x)$ of
index $\nu =(d - 2)/2$. When $d=1$,   $j_{-1/2}=\pi/2$.
\item[(ii)] 
Assume \eqref{eq:momentcondition},\eqref{eq:non-hom3}. If  $1/2<\rho<1$ ($ p_d<p<1$), then
\begin{equation}\label{eq:CLIL4}\liminf_{n\to \infty}\Big(\frac{\log\log n}{n}\Big)^{1/2} \max_{m\le n}\|\bm T_m-\mu_Z m^{\rho}\bm \xi\|=\frac{j_{\nu}}{\sqrt{2d}}(\ep[Z_1^2])^{1/2} \; a.s.
\end{equation}

\end{itemize}
\end{corollary}
  Interestingly, the limit values in the Chung type law of the iterated logarithm do not depend on the memory parameter $p$. In contrast, the limit values in the law of the iterated logarithm do. By taking $Z_i\equiv 1$, we obtain the Chung type law of the iterated logarithm of the elephant random walk $\bm S_n$. The Chung type law of the iterated logarithm of $\bm S_n$ is the same as that of a simple random walk for all $\rho<1/2$. 
 
 Corollary \ref{cor:CLIL} follows from the Gaussian approximation and the following small ball probability. 
 \begin{proposition}\label{prop:smallball1} Let $\bm W(t)$ and $\bm B(t)$ be two independent $d$-dimensional standard Brownian motions,
 \begin{equation}\label{eq:propsmallball1.1} \bm I(t)=\sigma_1 \int_0^t (t/s)^{\rho_1}d\bm W(t)+\sigma_2\int_0^t  (t/s)^{\rho_2}d\bm B(t), \;\; \rho_1,\rho_2<1/2. 
 \end{equation}
 Then 
  \begin{align}\label{eq:propsmallball1.2} 
   &\lim_{\epsilon\to 0}\epsilon^2 \log\pr\left(\sup_{0\le t\le 1}\|\bm I(t)\|<\epsilon\right)\nonumber\\
   = & 
   (\sigma_1^2+\sigma_2^2)\lim_{\epsilon\to 0}\epsilon^2 \log\pr\left(\sup_{0\le t\le 1}\|\bm B(t)\|<\epsilon\right)=-\frac{j_{(d-2)/2}^2}{2}(\sigma_1^2+\sigma_2^2). 
 \end{align}
 \end{proposition} 

\section{Applications to the  center of mass}
\setcounter{equation}{0}

In this section, we consider the  center of mass $\bm C_n$ of   $\bm T_n$  defined by
$$ \bm C_n =\frac{1}{n}\sum_{k=1}^n \bm T_k. $$
The question of the asymptotic behavior of the center of  a standard random walk in $\mathbb R^d$   was first raised by Paul Erd\"os.  
Lo and Wade \cite{LoWade2019} obtained the central limit theorem. Recently,  Bercu and  Laulin \cite{BercuLaulin2021} studied the asymptotic behavior of the center of mass of an ERW in $\mathbb R^d$ and obtained the rate of almost sure convergence, the central limit theorem and the law of the iterated logarithm. Also, they only obtained  an upper bound of the law of the iterated logarithm when $p<p_d$:
$$\limsup_{n\to \infty} \frac{\|\bm C_n\|}{\sqrt{2n\log\log n}}\le \frac{\sqrt{3}+\sqrt{1-2\rho}}{(1+\rho)\sqrt{3d(1-2\rho)}} \; a.s. $$

Now, by Theorem \ref{th:multi-RPWInva} and the properties of Gaussian processes, we have the following result on the asymptotic behavior of the center of mass $\bm C_n$ of an ERW with random step sizes in $\mathbb R^d$. By letting $Z_i=1$, then $\mu_Z=1$, $\sigma_Z=1$,   the central limit theorem and law of the iterated logarithm of the center of mass of $\bm S_n$  follow as a special case.  
 \begin{theorem}\label{cor:multi-CMERW} Consider the pair $\bm T_n$ and $\bm C_n$ of  an ERW with random step sizes    and its center of mass.  Let $\bm W(t)$,  $\bm B(t)$, $\bm G_t$ (when $\rho<1/2$) and $\widehat{\bm G}_t$ be defined as in Theorem \ref{th:multi-RPWInva}, and $\bm G_t$ is defined as in \eqref{eq:multi-thRPWInva15} when $\rho>1/2$.  Let
 \begin{equation}\label{eq:V-VmatrixofC}
 \bm \Lambda_{\rho}= \begin{cases} \begin{pmatrix}  \sigma_Z^2 +\frac{ \mu_Z^2 }{  1-2\rho } & 
 \frac{\sigma_Z^2}{2}+\frac{\mu_Z^2}{ (1-2\rho)(2-\rho)} \\
\frac{\sigma_Z^2}{2}+\frac{\mu_Z^2}{ (1-2\rho)(2-\rho)}  & \frac{\sigma_Z^2}{3}+\frac{2\mu_Z^2}{3(1-2\rho)(2-\rho)}\end{pmatrix} \otimes\frac{I_d}{d}, & \rho <1/2,\\
\frac{\mu_Z^2}{d}\begin{pmatrix} 1 & 2/3 \\ 2/3 & 4/9\end{pmatrix} \otimes I_d, & \rho=1/2, \\
\begin{pmatrix}  \sigma_Z^2 +\frac{ \mu_Z^2 }{  2\rho-1 } & 
 \frac{\sigma_Z^2}{2}+\frac{\mu_Z^2}{ (2\rho-1)(1+\rho)}  \\
\frac{\sigma_Z^2}{2}+\frac{\mu_Z^2}{ (2\rho-1)(1+\rho)}  & \frac{\sigma_Z^2}{3}+\frac{2\mu_Z^2}{3(2\rho-1)(1+\rho)}\end{pmatrix} \otimes\frac{I_d}{d}, & \rho>1/2.
\end{cases}
\end{equation}
\begin{itemize}
  \item[\rm (i)] 
Assume \eqref{eq:momentcondition} and \eqref{eq:non-hom}. If $\rho<1/2$ ($p<p_d$), then for some $\tau>0$, \eqref{eq:multi-thRPWInva2} holds and
\begin{align}\label{eq:multi-thCMERWInv1}
 \bm C_n=&\frac{1}{n}\left[\frac{\sigma_Z}{\sqrt{d}}\int_0^n  \bm W(t)dt+
\frac{\mu_Z}{ \sqrt{d }}   \int_0^n \bm G_t dt \right]
 +o\big(n^{1/2-\tau}\big) \; a.s.,
 \end{align}
\begin{align} \label{eq:multi-thCMERWCLT1}
&\frac{(\bm T_n, \bm C_n)}{\sqrt{n}}\overset{\mathscr{D}}\to N(\bm 0,\bm    \Lambda_{\rho}),
\end{align}
\begin{align} \label{eq:multi-thCMERWLIL1}
&\limsup_{n\to \infty} \frac{\|\bm C_n\|}{\sqrt{2n\log\log n}}= \sqrt{\frac{\sigma_Z^2}{3d}+\frac{2\mu_Z^2}{3d(1-2\rho)(2-\rho)}}\; a.s. 
\end{align}
  \item[\rm (ii)]  Assume $\ep Z_1^2<\infty$  and \eqref{eq:non-hom}.  If $\rho=1/2$ ($p=p_d$), then  
 \begin{equation} \label{eq:multi-thCMERWInv2} 
 \bm T_n=\frac{\mu_Z}{\sqrt{d}} \widehat{\bm G}_n+O(\sqrt{n\log\log n}) \text{ and } \bm C_n=\frac{2\mu_Z}{3\sqrt{d}} \widehat{\bm G}_n+O(\sqrt{n\log\log n}) \; a.s. 
 \end{equation}
 \begin{align} \label{eq:multi-thCMERWCLT2}
 \frac{(\bm T_n, \bm C_n)}{\sqrt{n\log n}} \overset{\mathscr{D}}\to N \Big(\bm 0,   \bm \Lambda_{1/2} \Big),
\end{align}
\begin{align} \label{eq:multi-thCMERWLIL2}
&\limsup_{n\to \infty} \frac{\|\bm C_n\|}{\sqrt{2n(\log n)(\log\log\log n)}}=  \frac{2|\mu_Z| }{3\sqrt{d}} \; a.s. 
\end{align}
\item[\rm (iii)] Assume $\ep Z_1^2<\infty$ and   \eqref{eq:non-hom2}. If $\rho>1/2$ ($p>p_d$), then  
  $$  \frac{\bm C_n}{n^{\rho}}\to \frac{\mu_Z\bm \xi}{1+\rho}\; a.s.  $$
  Furthermore, 
  $$  \left\|\frac{\bm C_n}{n^{\rho}}-    \frac{\mu_Z\bm \xi }{1+\rho} \right\|_q=O\big(n^{1/2-\rho}+|\sum_{i=n}^{\infty}(\rho_i-\rho)/i|\big), $$
  if $\ep|Z_1|^q<\infty$, $q\ge 2$, where $\bm\xi$ is the same as that in Theorem \ref{th:multi-RPWInva2}.
 \item[\rm (iv)] 
Assume \eqref{eq:momentcondition} and \eqref{eq:non-hom3}. If $1/2<\rho<1$ ($p_d<p<1$), then  for some $\tau>0$, \eqref{eq:multi-thRPWInva14} holds and
\begin{align}\label{eq:multi-thCMERWInv3}
 \bm C_n-  n^{\rho}\frac{\mu_Z\bm \xi}{1+\rho}=&\frac{1}{n}\left[\frac{\sigma_Z}{\sqrt{d}}\int_0^n  \bm W(t)dt+
\frac{\mu_Z}{ \sqrt{d }}   \int_0^n \bm G_t dt \right]
 +o\big(n^{1/2-\tau}\big) \; a.s., 
 \end{align}
\begin{align} \label{eq:multi-thCMERWCLT3}
&\frac{(\bm T_n-\mu_Z n^{\rho}\bm \xi, \bm C_n-\mu_Z n^{\rho}\bm \xi/(1+\rho))}{\sqrt{n}}\overset{\mathscr{D}}\to N(\bm 0,\bm    \Lambda_{\rho}),
\end{align}
\begin{align} \label{eq:multi-thCMERWLIL3}
&\limsup_{n\to \infty} \frac{\|\bm C_n-\mu_Z n^{\rho}\bm \xi/(1+\rho)\|}{\sqrt{2n\log\log n}}= \sqrt{\frac{\sigma_Z^2}{3d}+\frac{2\mu_Z^2}{3d(2\rho-1)(1+\rho)}}\; a.s. 
\end{align}
  \end{itemize}
\end{theorem}

Also,  in \eqref{eq:multi-thCMERWCLT1}, \eqref{eq:multi-thCMERWLIL1}, \eqref{eq:multi-thCMERWCLT3} and \eqref{eq:multi-thCMERWLIL3},  the condition $\ep|Z_1|^{2+\epsilon}<\infty$ can be removed. 

The following is the Chung type law of the iterated logarithm of   $\bm C_n$.  
\begin{corollary} \label{cor:CLIL2} \begin{itemize}
\item[(i)]Assume \eqref{eq:momentcondition} and \eqref{eq:non-hom}. If $\rho<1/2$ ($p<p_d$), then
\begin{equation}\label{eq:CLIL2.1}\liminf_{n\to \infty} \frac{(\log\log n)^{3/2}}{n^{1/2}} \max_{m\le n}\|\bm C_m\|=(3\kappa_d)^{3/2}(\ep Z_1^2/d)^{1/2} \;\; a.s., 
\end{equation}  
where  
\begin{equation}\label{eq:SBPofIB} \kappa_d=-\lim_{\epsilon\to 0}\epsilon^{2/3}\log \pr\left(\sup_{0\le t\le 1}\Big\|\int_0^t\bm B(s)ds\Big\|< \epsilon\right)  
\end{equation}
is a positive constant which only depends on $d$. 
\item[(ii)] 
Assume \eqref{eq:momentcondition} and \eqref{eq:non-hom3}. If  $1/2<\rho<1$ ($ p_d<p<1$), then
\begin{equation}\label{eq:CLIL2.2}\liminf_{n\to \infty} \frac{(\log\log n)^{3/2}}{n^{1/2}} \max_{m\le n}\left\|\bm C_m-m^{\rho}\frac{\mu_Z\bm\xi}{1+\rho}\right\|=(3\kappa_d)^{3/2}(\ep Z_1^2/d)^{1/2} \;\; a.s.
\end{equation} 
\end{itemize}
\end{corollary}

 The Chung type law of the iterated logarithm of the center of mass is new even for the simple random walk ($\rho=0$ and $Z_i\equiv 1$). Comparing Corollary \ref{cor:CLIL2} with Corollary \ref{cor:CLIL}, one can find that the center of mass slows down in the Chung type law of the iterated logarithm.

 Corollary \ref{cor:CLIL2}  follows from the Gaussian approximation and the following small ball probability and the Chung type law of the iterated logarithm  of the weighted maximum of integrated $\bm I_t$.  
 \begin{proposition}\label{prop:smallball2} Let $\bm W(t)$ and $\bm B(t)$ be two independent $d$-dimensional standard Brownian motions, $\bm I(t)$ is defined as in \eqref{eq:propsmallball1.1}, $\rho_1,\rho_2<1/2$.
 Then for $\alpha<3/2$, 
  \begin{align} 
  &\lim_{\epsilon\to 0}\epsilon^{2/3} \log\pr\left(\sup_{0< t\le 1}\Big\| \int_0^t\frac{1}{s^{\alpha}}\bm I(s)ds\Big\|<\epsilon\right)\nonumber\\
  =& \lim_{\epsilon\to 0}\epsilon^{2/3} \log\pr\left(\sup_{0< t\le 1}\Big\|\frac{1}{t^{\alpha}}\int_0^t\bm I(s)ds\Big\|<\epsilon\right) 
  =   -\frac{\kappa_d}{1-2\alpha/3}(\sigma_1^2+\sigma_2^2)^{1/3} \label{eq:propsmallball2.2} 
 \end{align}
 and
 \begin{align}\label{eq:propsmallball2.3} 
 &\liminf_{T\to \infty}\frac{(\log\log T)^{3/2}}{T^{3/2-\alpha}}\sup_{0< t\le T}\Big\| \int_0^t\frac{1}{s^{\alpha}}\bm I(s)ds\Big\|\nonumber\\
 =&\liminf_{T\to \infty}\frac{(\log\log T)^{3/2}}{T^{3/2-\alpha}}\sup_{0< t\le T}\Big\| \frac{1}{t^{\alpha}}\int_0^t\bm I(s)ds\Big\|
 =\Big(\frac{\kappa_d}{1-2\alpha/3}\Big)^{3/2}(\sigma_1^2+\sigma_2^2)^{1/2}\; \; a.s.
 \end{align}
 \end{proposition} 
\begin{remark} By \eqref{eq:propsmallball2.3}, the constant in Theorem 3 of Duker,  Li and  Linde \cite{DLL2000} is 
$C_{\alpha}= \big( \kappa_1/(1-2\alpha/3)\big)^{3/2}$. 
\end{remark}
 
 V\'azquez Guevara \cite{Vazque2019},  Bercu and  Laulin\cite{BercuLaulin2019,BercuLaulin2021} also obtained quadratic strong law for the ERW and the center of its mass, respectively. 
From  the Gaussian approximation,  we can obtain similar strong laws  for  various forms of functions much more general than  the quadratic function. The results are stated as   the following almost sure central limit theorem of   the ERW with random step sizes and the center of its mass. Here, we  use  a mild
regularity assumption of Ibragimov and Lifshits \cite{IbragimovLifshits1998} for the considered function. For the almost sure central limit theorems of i.i.d. random variables, we refer to Brosamler \cite{Brosamler1988}, Schatter \cite{Schatter1988}, Lacey and Philipp \cite{LaceyPhilipp1990},  Ibragimov and Lifshits \cite{IbragimovLifshits1998}, and  Berkes et al \cite{BCH1998}.

 \begin{corollary}\label{cor:ASCLT} Consider the pair $\bm T_n$ and $\bm C_n$ of  an ERW with random step sizes    and its center of mass. 
 Assume \eqref{eq:momentcondition} and \eqref{eq:non-hom}. Let $\bm \Lambda_{\rho}$  be defined as in \eqref{eq:V-VmatrixofC}. Let $A_0,H_0>0$, $\psi(x)$ be a continuous function such that $\psi(x)$ is  non-decreasing   on $[A_0,\infty)$ and $\psi(x)\exp\{-H_0x^2\}$ is non-increasing on $[A_0,\infty)$.  
  
\begin{itemize}
  \item[\rm (i)] If  for some $2d\times2d$   real matrix $\bm M$,
  \begin{equation}\label{eq:ASCLTcond1} \ep\Big[\psi\big(\|\widetilde{\bm C}\bm M\| \big)\Big]<\infty  \text{ with }  
 \widetilde{\bm C}\overset{\mathscr{D}} \sim   N(\bm 0,\bm \Lambda_{\rho}), 
 \end{equation}
 then for any almost everywhere continuous function $f(\bm x)$ which satisfies
 \begin{equation}\label{eq:ASCLTcond2}|f(\bm x)|\le \begin{cases} \psi(\|\bm x\bm M\|), & \|\bm x\bm M\|\ge A_0,  \\
 C_0, & \|\bm x\bm M\|< A_0,
 \end{cases}
 \end{equation}
   we have
\begin{align}\label{eq:ASCLT1.1}
 \frac{1}{\log n}\sum_{k=1}^n &\frac{f\big((\bm T_k,\bm C_k)/\sqrt{k}\big)}{k} \to \ep f(\widetilde{\bm C}) \;a.s. \;\; \text{ if }  \rho<1/2, \\
\label{eq:ASCLT1.2}
 \frac{1}{\log n}\sum_{k=1}^n & \frac{f\Big(\big(\bm T_k-k^{\rho}\mu_Z\bm\xi,\bm C_k-k^{\rho}\frac{\mu_Z\bm\xi}{1+\rho}\big)/\sqrt{k}\Big)}{k} \to \ep f(\widetilde{\bm C}) \;a.s. \\
  & \;\; \text{ if }  1/2<\rho<1 \text{ and } \eqref{eq:non-hom3}.\nonumber
\end{align} 
 Furthermore, if $f(\bm x)$ is an almost everywhere continuous function with $|f(\bm x)|\le C_0 e^{\gamma\|\bm x\bm M\|}$, $0<\gamma<\gamma_0$, where
\begin{equation}\label{eq:ASCLTcond3} \gamma_0=\sup\left\{\gamma>0: \ep\left[\exp\big\{\gamma\|\widetilde{\bm C}\bm M\|^2\big\}\right]<\infty\right\}, 
\end{equation}
  then the condition $\ep[|Z_1|^{2+\epsilon}]<\infty$ can be relaxed to $\ep[Z_1^2]<\infty$. 
  \item[\rm (ii)] If $\rho=1/2$ 
  and    for some $2d\times2d$   real matrix $\bm M$,
  \begin{equation}\label{eq:ASCLTcond3} \ep\Big[\psi\big(\|\widetilde{\bm C}\bm M\|\big)\Big]<\infty  \; \text{ with} 
\;\;\widetilde{\bm C}\overset{\mathscr{D}}=\frac{\mu_Z}{\sqrt{d}}(1,2/3)\otimes \bm B(1)\overset{\mathscr{D}}\sim  N \Big(\bm 0,   \bm \Lambda_{1/2}\Big), 
\end{equation}
 then for any   function $f(\bm x)$ which satisfies \eqref{eq:ASCLTcond2} and $\pr(\widetilde{\bm C}\in D_f)=0$, we have
  \begin{align}\label{eq:ASCLT2}
&\frac{1}{\log\log n}\sum_{k=2}^n \frac{f\big((\bm T_k,\bm C_k)/\sqrt{k\log k}\big)}{k\log k} \to \ep f(\widetilde{\bm C}) \;a.s., 
\end{align}
where $D_f$ is the set of discontinuities of $f$. Furthermore,  
the condition $\ep[|Z_1|^{2+\epsilon}]<\infty$ can be relaxed to $\ep[Z_1^2]<\infty$.
 \item[\rm (iii)] In (i) and (ii),  the exceptional set of probability $0$ can be chosen universally for all  $f$ which satisfy the given conditions. 
\end{itemize}
\end{corollary}
The almost sure central limit theorem for $(\bm T_n,\bm C_n)$ is equivalent to that \eqref{eq:ASCLT1.1},  \eqref{eq:ASCLT2} hold for every bounded continuous function $f$. 
By taking $f(\bm x, \bm y)=\|\bm y\|^2$,  $\bm y^{\prime}\bm y$, $\|\bm x\|^2$, $\bm x^{\prime}\bm x$ with $\bm x,\bm y\in\mathbb R^d$,  we obtain (2.3), (2.4), (2.10), (2.11) of   Bercu and   Laulin \cite{BercuLaulin2021}, and (3.2), (3.3), (3.7), (3.8) of Bercu and Laulin \cite{BercuLaulin2019}.  By taking $f(\bm x, \bm y)=\|\bm x\bm u^{\prime}\|^{2r}$ with $\bm x,\bm y\in\mathbb R^d$, we obtain (16), (18) of   V\'azquez Guevara \cite{Vazque2019}.

\section{Proofs}
\setcounter{equation}{0}
\subsection{Proof of the Gaussian approximation}

{\bf Proof of Theorem \ref{th:multi-RPWInva}.} Recall \eqref{eq:martingales.1}.  Notice that $\Delta \bm M_{n,1}$ is a sequence of bounded martingales. We have
$$ \bm M_{n,1}=\sum_{m=1}^n \Delta \bm M_{n,1}=O\big((n\log\log n)^{1/2}\big)\; a.s. \; \text{ and } =O(n^{1/2}) \text{ in } L_r $$
for all $r>0$. 
Let $\gamma_{m,n}=\prod_{i=m}^n(1+\rho_{i+1}/i)$, $\gamma_{m,m-1}=1$. Write $\log (1+x)=x-x^2h(x)$. Then $0<h(x)\le 2$ for $x\ge -1/2$.  Notice that   the condition \eqref{eq:non-hom} implies \eqref{eq:non-hom2}, and $|\rho_i|\le 1$. 
Thus 
$$\frac{\gamma_{2,n-1}}{n^\rho}=\exp\left\{\rho(\sum_{i=2}^{n-1}\frac{1}{i}-\ln n)+\sum_{i=2}^{n-1}\frac{\rho_{i+1}-\rho}{i}+\sum_{i=1}^{n-1}\frac{\rho_{i+1}^2}{i^2}h(\rho_{i+1}/i)\right\} \to C_0>0$$
and
$$
\frac{\gamma_{2,n-1}}{C_0n^\rho}=\exp\left\{\sum_{i=n}^{\infty} \frac{\rho_{i+1}-\rho}{i}+O(n^{-1}) \right\}
=1+O\Big(n^{-1}+|\sum_{i=n}^{\infty}(\rho_{i+1}-\rho)/i|\Big). $$
Hence $|\gamma_{m,n}|\le C (n/m)^{\rho}$.  
By \eqref{eq:SAforS1} and Lemma 3.2 of Bai, Hu and Zhang \cite{BHZ02},
\begin{align*}
\bm S_n=  & O\left(n^{\rho}\sum_{k=1}^n(k\log\log k)^{1/2}/(k^{\rho}k)\right)\;\; a.s.\\
=& \begin{cases} O\big( (n\log\log n)^{1/2}\big) \; a.s. & \rho<1/2\\
O\big( (n\log\log n)^{1/2}\log n\big) \; a.s. & \rho=1/2\\
O\big( n^{\rho} \big) \; a.s. & \rho>1/2
\end{cases}\\
=&o(n^{1-\epsilon}) \;\; a.s. \text{ for some } \epsilon>0 \text{ when } \rho<1.
\end{align*} 
Similarly,
\begin{align*}
\bm S_n=  & O_{L_r}\left(n^{\rho}\sum_{k=1}^nk^{1/2}/(k^{\rho}k)\right) 
=  \begin{cases} O( n^{1/2})   & \rho<1/2 \\
O( n^{1/2}\log n)   & \rho=1/2\\
O ( n^{\rho})   & \rho>1/2
\end{cases}\\
= &o(n^{1-\epsilon}) \;\;\text{ in } L_r \text{ for some } \epsilon>0 \text{ when } \rho<1.
\end{align*}  
For a vector $\bm v=(v_1,\ldots,v_d)$, we denote $\bm v^2=(v_1^2,\ldots,v_d^2)$. 
Similarly to \eqref{eq:condExp1}, we have that
\begin{align*}
 \ep[\bm\sigma_{n+1}^2|\beta_n]=&p_{n+1}\bm\sigma_{\beta_n}^2+\frac{1-p_{n+1}}{2d-1}\sum_{\bm\sigma\in \Sigma\setminus \bm\sigma_{\beta_n}}\bm\sigma^2  
 = p_{n+1}\bm\sigma_{\beta_n}+\frac{1-p_{n+1}}{2d-1}\left(\sum_{\bm\sigma\in \Sigma}\bm\sigma^2-\bm\sigma_{\beta_n}^2\right)\\
 =&p_{n+1}\bm\sigma_{\beta_n}^2+\frac{1-p_{n+1}}{2d-1}\left(2\bm 1-\bm\sigma_{\beta_n}^2\right)=\rho_{n+1} \bm\sigma_{\beta_n}^2+\frac{1-\rho_{n+1}}{d}\bm 1,
 \end{align*}
where $\bm 1=(1,\ldots,1)$. Thus, 
$$\ep[\bm\sigma_{n+1}^2|\mathscr{F}_n]=\rho_{n+1} \frac{\sum_{k=1}^n\bm \sigma_k^2}{n}+\frac{1-\rho_{n+1}}{d}\bm 1=\rho_{n+1} \frac{\sum_{k=1}^n\bm \sigma_k^2- \bm 1 n/d}{n}+\frac{1}{d}\bm 1. $$
Then, for $\bm q_n=\sum_{k=1}^n \bm\sigma_k^2-n\bm 1/d$ we also have the form of equation \eqref{eq:SAforS1}:
\begin{equation}\label{eq:var-sigma} \bm q_{n+1}=\bm q_n+\rho_{n+1} \frac{\bm q_n}{n}+\bm\sigma_{n+1}^2-\ep[\bm\sigma_{n+1}^2|\mathscr{F}_n].  
\end{equation}
Hence, we also have
$$ \sum_{k=1}^n \bm\sigma_k^2=n\frac{\bm 1}{d}+o(n^{1-\epsilon})\;\; a.s. \; \text{ and in } L_r \; \text{ when } \rho<1. $$

Now, we consider the variance-covariance matrix of the martingale differences $\Delta\bm M_n =(\Delta\bm M_{n,1},\Delta\bm M_{n,2})$. It is easily seen that
$\ep[\Delta\bm M_{n+1,1}^t\Delta\bm M_{n+1,2}|\mathscr{F}_n]=0$, and when $\rho<1$,
\begin{align*}
&\ep[\Delta\bm M_{n+1,2}^t\Delta\bm M_{n+1,2}|\mathscr{F}_n]
=\sigma_Z^2\ep[\bm\sigma_{n+1}^t\bm\sigma_{n+1}|\mathscr{F}_n]=\sigma_Z^2\ep[diag(\bm\sigma_{n+1}^2)|\mathscr{F}_n]\\
=&\sigma_Z^2 diag\left(\rho_{n+1}\frac{\sum_{k=1}^n\bm \sigma_k^2-n\bm 1/d}{n}+\frac{1}{d}\bm 1\right)
=\frac{\sigma_Z^2}{d}I_d+o(n^{-\epsilon})\; a.s. \; \text{ and in } L_1, 
\end{align*}
\begin{align}\label{eq:var-covar-M1}
 \ep[\Delta\bm M_{n+1,1}^t\Delta\bm M_{n+1,1}|\mathscr{F}_n]
= & \ep[diag(\bm\sigma_{n+1}^2)|\mathscr{F}_n]-\rho^2\big(  \frac{\bm S_n}{n}\big)^t\big(  \frac{\bm S_n}{n}\big)\nonumber \\
=&\frac{1}{d}I_d+o(n^{-\epsilon})\; a.s. \; \text{ and in } L_r,
\end{align}
 It follows that
\begin{equation}\label{eq:var-covar}
\ep[\Delta\bm M_{n+1}^t\Delta\bm M_{n+1}|\mathscr{F}_n]=\begin{pmatrix} 1 & 0\\
0 & \sigma_Z^2\end{pmatrix}\otimes \frac{I_d}{d}+o(n^{-\epsilon}) \; a.s. \; \text{ and in } L_1.
\end{equation}
The conditions in Theorems 3.1 and 2.2 of Zhang \cite{Zhang2016} are satisfied. 
\eqref{eq:multi-thRPWInva1}, \eqref{eq:multi-thRPWInva4} and \eqref{eq:multi-thRPWInva7} can follow from these two theorems. 
Here, we would   give a direct proof rather than  
derive (i)-(iii) from the complex formula of the limit process  for a general  recursive stochastic algorithm.  

Notice \eqref{eq:var-covar} and $\sup_n\ep[\|\Delta\bm M_n\|^{2+\epsilon}<\infty$, 
By Theorem 1.3 of Zhang
\cite{Zhang2004} or Theorem 1 of Eberlein \cite{Eberlein1986}, the $2d$ independent standard Brownian monitions $\{ B_k(t), W_k(t); t\ge 0, k=1,\ldots,d\}$ can be structured such that
\begin{equation}\label{eq:app-martingale} \bm M_{n,1}=\frac{1}{\sqrt{d}}\bm B(n)=o(n^{1/2-\tau})\;  a.s. \text{ and }\; 
\bm M_{n,2}=\frac{\sigma_Z}{\sqrt{d}}\bm W(n)=o(n^{1/2-\tau})\;  a.s. 
\end{equation}
for some $0<\tau<1$. Let $\widehat{\bm G}_t$ be defined as in \eqref{eq:multi-thRPWInva4}. Then
$$ \widehat{\bm G}_t=\rho\int_1^t \frac{\widehat{\bm G}_s}{s}ds +\bm B(t)-\bm B(1), $$
\begin{align}
 \widehat{\bm G}_{n}=&\rho\sum_{k=1}^{n-1}\frac{\widehat{\bm G}_k}{k}+\bm B(n)+o(n^{\delta})\; a.s., \forall \delta>0
 \nonumber\\
 =&  \sum_{k=1}^{n-1}\rho_{k+1}\frac{\widehat{\bm G}_k}{k}+\bm B(n)+o(n^{\delta})+\sum_{k=1}^{n-1}\frac{\rho_{k+1}-\rho }{k}\widehat{\bm G}_k. 
 \end{align}
Thus, by the above equality and  \eqref{eq:SAforS1},  $\bm S_{n}-\frac{1}{\sqrt{d}}\widehat{\bm G}_{n}$ also has the form of equation \eqref{eq:SAforS1}:
\begin{align*}
 &\bm S_{n}-\frac{1}{\sqrt{d}}\widehat{\bm G}_{n}= \sum_{k=1}^{n-1}\rho_{k+1}\frac{\bm S_k-\frac{1}{\sqrt{d}}\widehat{\bm G}_k}{k}+\bm\delta_{n,1}+\bm\delta_{n,2}, \\
 & \; \bm\delta_{n,2}= \sum_{k=1}^{n-1}\frac{\rho_{k+1}-\rho }{k}\widehat{\bm G}_k, \;\; \bm\delta_{n,1}=o(n^{1/2-\tau})  \;\; a.s. 
 \end{align*}
 We will show that under the condition \eqref{eq:non-hom}, $\bm \delta_{n,2}=o(n^{1/2-\epsilon})$ a.s. (c.f. Lemma \ref{lem1} below). 
 By Lemma 3.2 of Bai, Hu and Zhang \cite{BHZ02} again, we have
 $$ \bm S_{n}-\frac{1}{\sqrt{d}}\widehat{\bm G}_{n}=O\Big(n^{\rho}\sum_{k=1}^n\|\bm\delta_{k,1}+\bm\delta_{k,2}\|/(k^{\rho}k)\Big)=
 \begin{cases} o(n^{(1/2-\tau/2)\vee \rho}) \; a.s. & \text{ if } \rho<1/2,\\
 O(n^{\rho})\; a.s. &\text{ if } \rho=1/2.
 \end{cases} $$
\eqref{eq:multi-thRPWInva1},\eqref{eq:multi-thRPWInva2}, \eqref{eq:multi-thRPWInva4}, \eqref{eq:multi-thRPWInva5} are proven if noting that $\int_0^1s^{-\rho}d\bm B(s)$ is finite when $\rho<1/2$.

(iii) is implied by Theorem \ref{th:multi-RPWInva2}. $\Box$

{\bf Proof of Theorem \ref{th:multi-RPWInva2}.} 
Recall $\gamma_{m,n}=\prod_{i=m}^n(1+\rho_{i+1}/i)$, $\gamma_{m,m-1}=1$. Then
\begin{align}\label{eq:Expan-S}
\bm S_n=& \bm S_{n-1}\gamma_{n-1,n-1}+\Delta \bm M_{n,1}\nonumber\\
=&\bm S_{n-2}\gamma_{n-2,n-1}+\Delta\bm M_{n-1,1}\gamma_{n-1,n-1}+
\Delta \bm M_{n,1}\gamma_{n,n-1} \nonumber\\
=&\bm S_1\gamma_{1,n-1}+\sum_{k=2}^n\Delta\bm M_{k,1} \gamma_{k,n-1}\nonumber\\
=&\big(\bm S_1(1+\rho_2)   +\sum_{k=2}^n\Delta\bm M_{k,1} \gamma_{2,k-1}^{-1}\big)\gamma_{2,n-1}.
\end{align}

We prove the theorem via four steps. 

{\em Step 1}. We show (a).
For the martingale differences $\{ \Delta\bm M_{k,1}/ \gamma_{2,k-1};k\ge 2\}$ we have
$$ \sum_{k=2}^{\infty}\frac{\ep\big[\|\Delta\bm M_{k,1}\|^2\big|\mathscr{F}_{k-1}\big]}{ \gamma_{2,k-1}^2}
\le c\sum_{k=1}^{\infty} k^{-2\rho}<\infty. $$
It follows that (c.f. Theorem 2.15 of Hall and Heyde \cite{HallHeyde})
$$ \bm\eta=\bm S_1(1+\rho_2) +\sum_{k=1}^{\infty}\frac{\Delta\bm M_{k,1}}{ \gamma_{2,k-1}} \text{ is a.s. finite}. $$
 Thus
$$ \frac{\bm S_n}{n^{\rho}}=\frac{\bm S_n}{\gamma_{2,n-1}}\frac{\gamma_{2,n-1}}{n^{\rho}}\to  C_0\bm \eta=:\bm \xi  \;\; a.s. $$
Write $b_n^2=\sum_{k=n+1}^{\infty}\gamma_{2,k-1}^{-2}$,
$$\bm\eta_n=\frac{\bm S_n}{\gamma_{2,n-1}}=\bm S_1(1+\rho_2)  +\sum_{k=2}^n\frac{\Delta\bm M_{k,1}}{\gamma_{2,k-1}}. $$
Then $b_n^2\asymp n^{1-2\rho}$,
$$ \ep[\|\bm\eta-\bm\eta_n\|^r]\le C_r\ep\left(\sum_{k=n+1}^{\infty}\ep\big[\|\Delta\bm M_{k,1}\|^2\big|\mathscr{F}_{k-1}\big]\gamma_{2,k-1}^{-2}\right)^{r/2}
\le C_r b_n^r, $$
$$ \big\|\frac{\bm S_n}{n^{\rho}}-C_0\bm\eta_n\big\|_r=O\big(\big|\frac{\gamma_{2,n-1}}{n^{\rho}}-C_0\big|\Big)
=O\Big(\big|\sum_{i=n}^{\infty}(\rho_i-\rho)/i\big|\Big).$$
It follows that
$$\big\|\frac{\bm S_n}{n^{\rho}}-C_0\bm\eta\big\|_r=O\left(n^{1/2-\rho}+\big|\sum_{i=n}^{\infty}(\rho_i-\rho)/i\big|\right).$$
By noting that 
$$ \sum_{k=1}^n \bm\sigma_k(Z_k-\mu_Z)=O(\sqrt{n\log\log n}) \; a.s. $$
$$(\ep \|\sum_{k=1}^n \bm\sigma_k(Z_k-\mu_Z)\|^q)^{1/q}=O(\sqrt{n})\; \text{ if }\; \ep |Z_1|^q<\infty,\; q\ge 2. $$
\eqref{eq:multi-thRPWInva7} and \eqref{eq:multi-thRPWInva8} are proven. 

{\em Step 2}. We show (b).    By \eqref{eq:Expan-S} we have $\ep[\bm S_n]=\ep[\bm S_1]\gamma_{1,n-1}=0$. Thus,
$\ep[\bm\xi]=0$. Notice
\begin{align*}
 \ep\left[\bm S_{n+1}^{\prime} \bm S_{n+1}|\mathscr{F}_n\right]= &  \bm S_n^{\prime} \bm S_n 
+2 \bm S_n^{\prime}\ep\left[ \bm \sigma_{n+1}|\mathscr{F}_n\right]+\ep\left[diag(\bm\sigma_{n+1}^2)|\mathscr{F}_n\right]\\
=& \big(1+\frac{2\rho_{n+1}}{n}\big) \bm S_n^{\prime} \bm S_n 
+diag\left(\rho_{n+1}\frac{\sum_{k=1}^n\bm \sigma_k^2-n\bm 1/d}{n}+\frac{1}{d}\bm 1\right).
\end{align*}
By \eqref{eq:var-sigma}, 
$$ \ep\left[\sum_{k=1}^n \bm\sigma_k^2-n\bm 1/d\right]=\big(1+\frac{\rho_{n}}{n-1}\big)\cdots \big(1+\frac{\rho_2}{1}\big)\ep[\bm\sigma_1^2-\bm 1/d]=0. $$
It follows that
$$  \ep\left[\bm S_{n+1}^{\prime} \bm S_{n+1}\right]=\big(1+\frac{2\rho_{n+1}}{n}\big) \ep\left[\bm S_n^{\prime} \bm S_n\right]
+\frac{I_d}{d}. $$ 
Similar to \eqref{eq:Expan-S}, we have
\begin{align*}
\ep\left[\bm S_n^{\prime} \bm S_n\right]=&\ep\left[\bm S_1^{\prime} \bm S_1\right]\gamma_{1,n-1}^{\ast}+\frac{I_d}{d}\sum_{k=1}^n  \gamma_{k,n-1}^{\ast}\\
=& \frac{I_d}{d}\left[(1+2\rho_2)(1+\rho_3) + (1+\rho_3) +\sum_{k=3}^n \frac{1}{\gamma_{3,k-1}^{\ast}}\right]\gamma_{3,n-1}^{\ast}.
\end{align*}
where $\gamma_{m,n}^{\ast}=\prod_{i=m}^n(1+2\rho_{i+1}/i)$, $\gamma_{m,m-1}^{\ast}=1$, $\gamma_{3,n-1}^{\ast}/n^{2\rho}\to C^{\ast}$. Thus
\begin{align*}
\frac{\ep\left[\bm S_n^{\prime} \bm S_n\right]}{n^{2\rho}}
\to \frac{1}{d} C^{\ast}\left[2(2+\rho_2)(1+\rho_3)+\sum_{k=3}^{\infty}\frac{1}{\gamma_{3,k-1}^{\ast}}\right] I_d:=CI_d.
\end{align*}
\eqref{eq:multi-thRPWInva12} is proven.  

{\em Step 3}. We Show (c).   For the martingale
$\bm\eta-\bm\eta_n$,   \eqref{eq:var-covar} remains true since $\rho<1$. Thus, 
$$  b_n^{-2}
\sum_{k=n+1}^{\infty}\ep\big[(\Delta\bm M_{k,1})^{\prime}\Delta\bm M_{k,1}\big|\mathscr{F}_{k-1}\big]\gamma_{2,k-1}^{-2}
\to \frac{I_d}{d} \; a.s., $$
$$  b_n^{-3}
\sum_{k=n+1}^{\infty}\ep\big[\|\Delta\bm M_{k,1}\|^3\big|\mathscr{F}_{k-1}\big]\gamma_{2,k-1}^{-3}
\le Cn^{-3(1-2\rho)/2}\sum_{k=n+1}^{\infty} k^{-3\rho}\to 0. $$
The conditional Lindeberg condition is satisfied. By the conditional central limit theorem for martingales,   given $\mathscr{F}_n$, the conditional distribution of 
$b_n^{-1}(\bm \eta-\bm\eta_n)$ 
will almost surely converge to a multi-normal distribution $N(\bm 0, I_d/d)$. It follows that
$$ \ep\left[\exp\{ {\rm i}b_n^{-1}(\bm \eta-\bm\eta_n)\bm t^{\prime}\}\big|\mathscr{F}_n\right]\to e^{-\frac{\|\bm t\|^2}{2d}}\;\; a.s. $$
When $\bm t\ne 0$, let $E=I\{\bm \eta\bm t^{\prime}=x\}$, $I_n=\ep[I_E|\mathscr{F}_n]$. Then $I_n\to I_E$ a.s. and $\ep[|I_n-I_E||\mathscr{F}_n]\to 0$ in $L_1$. 
Thus,
\begin{align*}
&\lim_{n\to \infty} \ep\left[\exp\{ {\rm i}b_n^{-1}(\bm \eta-\bm\eta_n)\bm t^{\prime}\}I_E\big|\mathscr{F}_n\right]
=\lim_{n\to \infty} \ep\left[\exp\{ {\rm i}b_n^{-1}(\bm \eta-\bm\eta_n)\bm t^{\prime}\}I_n\big|\mathscr{F}_n\right]\\
& \;\; =  \lim_{n\to \infty} \ep\left[\exp\{ {\rm i}b_n^{-1}(\bm \eta-\bm\xi_n)\bm t^{\prime}\}\big|\mathscr{F}_n\right]I_n
=e^{-\frac{\|\bm t\|^2}{2d}}I_E \; a.s. 
\end{align*}
On the other hand, on the event $E$, $b_n^{-1}(\bm \eta-\bm\eta_n)\bm t^{\prime}=b_n^{-1}(x-\bm\eta_n \bm t^{\prime})$ is a $\mathscr{F}_n$-measurable random variables. By noticing $|e^{{\rm i}y}|=1$, it follows that
\begin{align*}
I_E=&\lim_{n\to \infty}\ep[I_E|\mathscr{F}_n]=
 \lim_{n\to \infty} \left|\exp\{ {\rm i}b_n^{-1}(\bm \eta-\bm\eta_n)\bm t^{\prime}\}\ep\left[I_E\big|\mathscr{F}_n\right]\right|\\
=&\lim_{n\to \infty} \ep\left[\exp\{ {\rm i}b_n^{-1}(\bm \eta-\bm\eta_n)\bm t^{\prime}\}I_E\big|\mathscr{F}_n\right]
=e^{-\frac{\|\bm t\|^2}{2d}}I_E \; a.s., 
\end{align*}
which implies that $I_E=0$ a.s. That is, $\pr(\bm\eta\bm t^{\prime}=x)=0$. Thus, $\pr(\bm \xi\bm t^{\prime}=x)=0$ for all $x$.

{\em Step 4}. We show (d).  When \eqref{eq:non-hom3} is satisfies, then 
$$ \sum_{i=n}^{\infty}\frac{(\rho_{i+1}-\rho)}{i^{1+\gamma}}=O(n^{-\epsilon_0-\gamma}), \;\; \gamma\ge 0. $$
In fact, let $T_i=\sum_{j=1}^i(\rho_{j+1}-\rho)$.
Then
\begin{align*}
\sum_{i=n}^{\infty}\frac{(\rho_{i+1}-\rho)}{i^{1+\gamma}}=& \sum_{i=n}^{\infty}T_i
\left(\frac{1}{i^{1+\gamma}}-\frac{1}{(i+1)^{1+\gamma}}\right)-\frac{T_{n-1}}{n^{1+\gamma}}\\
=& \sum_{i=n}^{\infty}\frac{O(i^{1-\epsilon_0})}{ii^{1+\gamma}}+\frac{O(n^{1-\epsilon_0})}{n^{1+\gamma}}=O(n^{-\epsilon_0-\gamma}).
\end{align*}
Thus, $\gamma_{2,n-1}=C_0n^{\rho}\big(1+O(n^{-\epsilon_0}))$. By noting that  $\bm M_{k,1}=\frac{1}{\sqrt{d}}\bm B(k)+o(k^{1/2-\tau})$ a.s., we have
\begin{align*}
& C_0(\bm\eta-\bm\eta_n)= \sum_{k=n+1}^{\infty}\frac{C_0\Delta\bm M_{k,1}}{ \gamma_{2,k-1}}
=\sum_{k=n+1}^{\infty}\frac{C_0\bm M_{k,1}-C_0\bm M_{k-1,1}}{\gamma_{2,k-1}}\\
=&\sum_{k=n+1}^{\infty}\Big(\frac{1}{ \gamma_{2,k-1}}-\frac{1}{ \gamma_{2,k}}\Big)C_0\bm M_{k,1}-\frac{C_0\bm M_{n,1}}{\gamma_{2,n}} \\
=& \frac{1}{\sqrt{d}}\sum_{k=n+1}^{\infty} \frac{\rho_{k+1} C_0 \bm B(k)+o(k^{1/2-\tau})}{k \gamma_{2,k}} -\frac{1}{\sqrt{d}}\frac{C_0\bm B(n)+o(n^{1/2-\tau})}{\gamma_{2,n}} \\
=& \frac{1}{\sqrt{d}}\sum_{k=n+1}^{\infty} \frac{\rho   \bm B(k)+o(k^{1/2-\tau})+(\rho_{k+1}-\rho)\bm B(k)}{k k^{\rho}(1+O(k^{-\epsilon_0}))} -\frac{1}{\sqrt{d}}\frac{\bm B(n)+o(n^{1/2-\tau})}{n^{\rho}(1+O(k^{-\epsilon_0}))} \\
=&\frac{1}{\sqrt{d}} \rho\int_n^{\infty} s^{-1-\rho} \bm B(s)ds-\frac{1}{\sqrt{d}}\frac{\bm B(n)}{n^{\rho}}  +\sum_{k=n+1}^{\infty} \frac{(\rho_{k+1}-\rho)\bm B(k)}{k k^{\rho}} +o(n^{1/2-\rho-\tau}) \\
=& \frac{1}{\sqrt{d}}\int_n^{\infty} s^{-\rho} d\bm B(s) + \bm\zeta_n +o(n^{1/2-\rho-\tau}) \; a.s.,
\end{align*}
where
\begin{align*}
\bm\zeta_n=& \sum_{k=n+1}^{\infty} \frac{(\rho_{k+1}-\rho)\bm B(k)}{ k^{1+\rho}}\\
=& \bm B(n)\sum_{k=n+1}^{\infty}\frac{\rho_{k+1}-\rho}{ k^{1+\rho}}+\sum_{j=n+1}^{\infty} \big(\bm B(j)-\bm B(j-1)\big) \sum_{k=j}^{\infty}\frac{\rho_{k+1}-\rho}{ j^{1+\rho}} 
\end{align*}
is a mean zero Gaussian random variable with
\begin{align*}
\ep\left[\|\bm\zeta_n\|^2\right]=d n \left(\sum_{k=n+1}^{\infty}\frac{\rho_{k+1}-\rho}{ k^{1+\rho}}\right)^2
+d \sum_{j=n+1}^{\infty} \left(\sum_{k=j}^{\infty}\frac{\rho_{k+1}-\rho}{ j^{1+\rho}} \right)^2=O(n^{1-2\rho-2\epsilon_0}).
\end{align*} 
It follows that
\begin{align*}
\sum_{n=1}^{\infty} \pr\left(\|\bm\zeta_n\|\ge \epsilon n^{1/2-\rho-\tau}\right)
\le\sum_{n=1}^{\infty}  \exp\Big\{- \frac{\epsilon^2n^{1 -2\rho-2\tau}}{2O(n^{1-2\rho-2\epsilon_0})}\Big\}<\infty,
\end{align*}
which implies $\bm\zeta_n=o(n^{1/2-\rho-\tau})$ a.s.
Thus,
\begin{align*}
  \frac{\bm S_n}{n^{\rho}}-  \bm \xi = & C_0(\bm\eta_n-  \bm \eta)+\left(\frac{\gamma_{2,n-1}}{n^{\rho}}-C_0\right) \bm \eta_n\\
=& -\frac{1}{\sqrt{d}}\int_n^{\infty} s^{-\rho} d\bm B(s) +o(n^{1/2-\rho-\tau}) +O(n^{-\epsilon_0})\; a.s.
\end{align*}
\eqref{eq:multi-thRPWInva13} holds. \eqref{eq:multi-thRPWInva14} follows from \eqref{eq:multi-thRPWInva13} and \eqref{eq:app-martingale}.
the proof is now completed.  $\Box$

Filially, we show a lemma which is used in the proof of Theorem \ref{th:multi-RPWInva}.
\begin{lemma}\label{lem1} Under the condition \eqref{eq:non-hom},
$$\bm\delta_{n,2}= \sum_{k=1}^{n-1}\frac{\rho_{k+1}-\rho }{k}\widehat{\bm G}_k=o(n^{1/2-\epsilon}) \;\; a.s.\; \text{ and in } L_2
\; \text{ for some } \epsilon>0.
$$
\end{lemma}
{\bf Proof.} Define $\rho_x=\rho_{k+1}$ for $x\in (k,k+1]$,
$$ \beta_{\rho}(s)=\int_s^{\infty}\frac{\rho_x-\rho}{x^{1-\rho}} dx, \;\rho\le 0 \text{ and } \beta_{\rho}(s)=-\int_1^s\frac{\rho_x-\rho}{x^{1-\rho}} dx, \;\rho> 0.   $$
Then
\begin{align*}
\bm\delta_{n,2}=&\sum_{k=1}^{n-1}\frac{\rho_{k+1}-\rho }{k^{1-\rho}} \int_1^k s^{-\rho} d\bm B(s)\\
=& \int_1^n \left[\frac{\rho_x-\rho}{x^{1-\rho}} \int_1^x s^{-\rho} d\bm B(s)\right]dx  +o(n^{\delta}) \; a.s. \; \text{ and in } L_2 \;\forall \delta>0\\
=&\int_1^ns^{-\rho}d\bm B(s)\int_s^n \frac{\rho_x-\rho}{x^{1-\rho}}dx  +o(n^{\delta})\\
=& \int_1^n \beta_{\rho}(s)s^{-\rho} d\bm B(s)-\beta_{\rho}(n)\int_1^n s^{-\rho}d\bm B(s)+o(n^{\delta})  
\end{align*}
Notice that $\{\int_1^tf(s)d\bm B(s);t\ge 1\}\overset{\mathscr{D}}=\{\bm B \big(\int_1^tf^2(s)ds\big); t\ge 1\}$. By the law of iterated logarithm of the Borwnian motion. It is sufficient to show that
 $ \beta_{\rho}(t)=o(t^{\rho-\epsilon})\; \text{ as } t\to \infty. $ 
First, $\alpha(t)=:\int_0^t(\rho_x-\rho)=o(t^{1-\epsilon})$. Then
\begin{align*}
 \beta_{0}(t)=&\int_t^{\infty}x^{-1}d\alpha(x)=   \alpha(x)x^{-1}\big|_{t}^{\infty}+\int_t^{\infty}\alpha(x)x^{-2}dx \\
 =& o(t^{-\epsilon})+\int_t^{\infty} o(x^{-\epsilon-1})dx=o(t^{-\epsilon}).
 \end{align*}
  When $\rho>0$, without loss of generality, we assume $\epsilon<\rho$.  Then, when $\rho< 0$, 
\begin{align*}
 \beta_{\rho}(t)=&-\int_t^{\infty}x^{\rho}d\beta_{0}(t)=   -\beta_0(x)x^{\rho}\big|_{t}^{\infty}+\rho\int_t^{\infty}\beta_0(x)x^{\rho-1}dx \\
 =& o(t^{\rho-\epsilon})+\int_t^{\infty} o(x^{\rho-\epsilon-1})dx=o(t^{\rho-\epsilon}),
 \end{align*}
 and,  when $\rho> 0$, 
\begin{align*}
 \beta_{\rho}(t)=&\int_1^tx^{\rho}d\beta_0(x)=   \beta_0(x)x^{\rho}\big|_{t}^{\infty}
 -\rho\int_1^t\beta_0(x)x^{\rho-1}dx \\
 =& o(t^{\rho-\epsilon})-\rho\int_1^t o(x^{\rho-\epsilon-1})dx=o(t^{\rho-\epsilon}).
\end{align*}
The proof is completed. $\Box$

\begin{remark} Under the condition \eqref{eq:non-hom2}, we have $\beta_0(t)\to 0$ as $\to \infty$. Then $\beta_{\rho}(t)=o(t^{\rho})$. Thus, for $\rho<1/2$,
$$ \bm\delta_{n,2}=o(\sqrt{n\log\log n})\;\; a.s.\; \text{ and }=o(\sqrt{n}) \text{ in } L_2, $$
which will imply
$$ \bm S_{n}-\frac{1}{\sqrt{d}} \bm G_{n}=o(\sqrt{n\log\log n})\;\; a.s.\; \text{ and }=o(\sqrt{n}) \text{ in } L_2. $$
Thus, for \eqref{eq:multi-thERWCLT1} and \eqref{eq:multi-thERWLIL1}, the condition \eqref{eq:non-hom} can be weakened to \eqref{eq:non-hom2}.
\end{remark}

\subsection{Proof of the CLT and LIL}

{\bf Proof of Corollary \ref{cor:multi-ERW}.} For \eqref{eq:multi-thERWCLT1} and \eqref{eq:multi-thERWCLT3}, it is sufficient to notice that
$$ \Var\{\bm G_t\}=I_d\Var\{t^{\rho}\int_0^ts^{-\rho}d B_1(t)\}=I_d t^{2\rho}\int_0^1 s^{-2\rho}ds=\frac{t}{1-2\rho}I_d, \; \rho<1/2$$
and
$$ \Var\{\bm G_t\}=I_d\Var\{\int_t^{\infty}(t/s)^{\rho}d B_1(s)\}=I_dt^{2\rho}\int_t^{\infty} s^{-2\rho}ds=\frac{t}{2\rho-1}I_d, \; \rho>1/2. $$
For \eqref{eq:multi-thERWCLT2} and \eqref{eq:multi-thERWLIL2}, it is sufficient to notice that
$$  \bm W(n)=O(\sqrt{n\log\log n}) \; a.s. $$
and 
\begin{equation}\label{eq:cormulti-ERW.1}\{\widehat{\bm G}_t;t\ge 1\}\overset{\mathscr{D}}=\{t^{1/2}\bm B(\log t); t\ge 1\} \text{ when } \rho=1/2.
\end{equation}
For \eqref{eq:multi-thERWLIL1} and \eqref{eq:multi-thERWLIL3}, it is sufficient to consider the   law of the iterated logarithm of the Gaussian processes on the right hands of \eqref{eq:multi-thRPWInva1} and \eqref{eq:multi-thRPWInva14}.   Proposition 5.10 of Zhang \cite{Zhang2012} gave the law of the iterated logarithm of similar kinds of such Gaussian processes.   We omit the details here  and pay the main attention on the same properties of the center of mass $\bm C_n$.   $\Box$

\bigskip
{\bf Proof of Theorem \ref{cor:multi-CMERW}.} We fist consider (ii). In this case, $\rho=1/2$, and we will still have \eqref{eq:multi-thRPWInva1} in spite  of that \eqref{eq:momentcondition} is not satisfied.   Thus,
$$ \bm T_n=O((n\log\log n)^{1/2})+\mu_Z\bm S_n=\frac{\mu_Z}{\sqrt{d}} \widehat{\bm G}_n+O((n\log\log n)^{1/2})\; a.s. $$
and
\begin{align*}
\bm C_n=&\frac{\mu_Z}{\sqrt{d}}\frac{1}{n}\sum_{k=1}^n \widehat{\bm G}_k+O(\sqrt{n\log\log n}) \; a.s.\\
=& \frac{\mu_Z}{\sqrt{d}}\frac{1}{n}\int_1^n \widehat{\bm G}_tdt+O(\sqrt{n\log\log n})\; a.s. \\
=& \frac{2\mu_Z}{3\sqrt{d}}\left(\widehat{\bm G}_n -\frac{1}{n}\int_1^ns d\bm B(s)\right)+O(\sqrt{n\log\log n}) \; a.s.
\end{align*}
It easily seen that 
\begin{equation}\label{eq:proof-CMERW1}\left\{\int_0^ts^{-\gamma} d\bm B(s);t\ge 0\right\}\overset{\mathscr{D}}=\left\{\bm B(t^{1-2\gamma}/(1-2\gamma)); t\ge 0\right\}, \;\gamma<1/2,
\end{equation}
and thus 
\begin{equation}\label{eq:proof-CMERW2}\int_1^ns^{-\gamma} d\bm B(s)=O(\sqrt{n^{1-2\gamma}\log\log n})\;\; a.s.
\end{equation}
In particular, $\int_1^nsd \bm B(s)=O(\sqrt{n^3\log\log n})$ a.s. It follows that
$$ \bm C_n=\frac{2\mu_Z}{3\sqrt{d}} \widehat{\bm G}_n+O(\sqrt{n\log\log n}) \; a.s. $$
\eqref{eq:multi-thCMERWInv2} is proven, and \eqref{eq:multi-thCMERWCLT2} and \eqref{eq:multi-thCMERWLIL2} follow immedaitely by \eqref{eq:cormulti-ERW.1}.
 
Next we consider (i). In this case,  $\rho<1/2$ case. By \eqref{eq:multi-thRPWInva5}, 
\begin{align*}
 \bm C_n=&\frac{1}{n}\sum_{k=1}^n \left[\frac{\sigma_Z}{\sqrt{d}}\bm W(k)+\frac{\mu_Z}{\sqrt{d }}  \bm G_k\right]+o(n^{1/2-\tau}) \; a.s.\\
 =&\frac{1}{n}  \left[\frac{\sigma_Z}{\sqrt{d}}\int_0^n \bm W(s)ds+\frac{\mu_Z}{\sqrt{d }}\int_0^n  \bm G_sds\right]+o(n^{1/2-\tau}) \; a.s. 
\end{align*}
\eqref{eq:multi-thCMERWInv1} is proved. Next, it is sufficient to derive the variance-covariance matrix  of the Gaussian process
$(\widetilde{\bm T}_T/T^{1/2}, \widetilde{\bm C}_T/T^{3/2})$ and the law of the iterated logarithm of $\widetilde{\bm C}_T$, where
\begin{equation}\label{eq:CMEWGauss1}
\widetilde{\bm T}_T=   \sigma_1 \bm W(T)+ \sigma_2 \bm G_T \;\; \text{ and }\;\;
\widetilde{\bm C}_T=   \sigma_1 \int_0^T \bm W(s)ds+\sigma_2\int_0^T  \bm G_sds. 
\end{equation}
 Note 
$$  \int_0^T\left[t^{\gamma}\int_0^ts^{-\gamma}d\bm B(s)\right]dt=
\begin{cases} \frac{1}{1+\gamma}\int_0^T\big(T^{1+\gamma}s^{-\gamma}-s\big) d\bm B(s), &\gamma<1/2,\gamma\ne -1\\
\int_0^T(\log T-\log s)s d\bm B(s), &\gamma=-1
\end{cases}. $$
Next, we assume $\rho\ne -1$ since the proof for the case $\rho=-1$ is the same. Then
\begin{align}\label{eq:CMEWGauss2}
\widetilde{\bm T}_T 
= & \sigma_1\int_0^T d \bm W(s)+ \sigma_2 \int_0^T  T^{\rho}s^{-\rho} d\bm B(s), \nonumber\\
\widetilde{\bm C}_T 
= & \sigma_1\int_0^T(T-s)d \bm W(s)+\frac{\sigma_2}{1+\rho}\int_0^T \big(T^{1+\rho}s^{-\rho}-s\big)d \bm B(s) 
\end{align} with
\begin{align}
\Var\{\widetilde{\bm T}_T/T^{1/2}\}=&\frac{1}{T}\left[\sigma_1^2\int_0^Tds+ \sigma_2^2 \int_0^T T^{2\rho}s^{-2\rho}ds\right]I_d \nonumber\\
=&\Big(\sigma_1^2  +\frac{ \sigma_2^2}{ 1-2\rho }\Big)I_d,
\label{eq:varinceGauss1}\\
\Cov\{\widetilde{\bm T}_T/T^{1/2},\widetilde{\bm C}_T/T^{3/2}\}=&\frac{1}{T^2}\left[\sigma_1^2\int_0^T(T-s)ds+ \frac{\sigma_2^2}{1+\rho}\int_0^T\big(T^{1+\rho}s^{-\rho}-s\big)T^{\rho}s^{-\rho}ds\right]I_d\nonumber\\
=&  \left[\frac{\sigma_1^2}{2} +\frac{\sigma_2^2}{ (1-2\rho)(2-\rho)} \right]I_d,
\label{eq:varinceGauss2}\\
\Var\{\widetilde{\bm C}_T/T^{3/2}\}=&\frac{1}{T^3}\left[\sigma_1^2\int_0^T(T-s)^2ds+ \frac{\sigma_2^2}{(1+\rho)^2}\int_0^T\big(T^{1+\rho}s^{-\rho}-s\big)^2ds\right]I_d \nonumber\\
=&  \left[\frac{\sigma_1^2}{3} +\frac{2\sigma_2^2}{3(1-2\rho)(2-\rho)} \right]I_d=:\sigma^2_{\rho} I_d.
\label{eq:varinceGauss3}
\end{align}
\eqref{eq:multi-thCMERWCLT1} is proven.  

For \eqref{eq:multi-thCMERWLIL1}, it is sufficient to show that
\begin{equation}\label{eq:LILforC}
\limsup_{T\to \infty} \frac{\|\widetilde{\bm C}_T\|}{T^{3/2}\sqrt{2\log\log T}}=\sigma_{\rho}\; a.s. 
\end{equation}
Note $\|\widetilde{\bm C}_s\|=\sup\{\widetilde{\bm C}_s\bm b^{\prime}:\|\bm b\|=1\}$. 
Denote $\|\widetilde{\bm C}_T\|_T=\sup_{0\le s\le T}\|\widetilde{\bm C}_s\|=\sup\{\widetilde{\bm C}_s\bm b^{\prime}:\|\bm b\|=1,0\le t\le T\}$. By the Borell inequality (c.f. Borell\cite{Borell1975}, Pisier\cite{Pisier1986}),
$$ \pr\left(\|\widetilde{\bm C}\|_T-\ep[\|\widetilde{\bm C}\|_T]\ge x\right)\le \exp\left\{-\frac{x^2}{2\sigma_T^2}\right\}, $$
where
\begin{align*}
 \sigma_T^2=&\sup\{\Var(\widetilde{\bm C}_s\bm b^{\prime}):\|\bm b\|=1,0\le t\le T\}\\
 =& \sup\{s^3\sigma^2_{\rho}\|\bm b\|^2:\|\bm b\|=1,0\le t\le T\}=T^3\sigma^2_{\rho}.
 \end{align*}
On the other hand, 
$$ \ep[\|\widetilde{\bm C}\|_T]\le |\sigma_1|\int_0^T\ep[\|\bm W(t)\|]dt+|\sigma_2|\int_0^T\ep[\|\bm G_t\|]dt
\le c\int_0^T\sqrt{t}dt \le c_0T^{3/2}. $$
Thus,
$$\pr\left(\|\widetilde{\bm C}\|_T\ge x+c_0 T^{3/2}\right)\le \exp\left\{-\frac{x^2}{2T^{3/2}\sigma^2_{\rho}} \right\}. $$
Choose $T_k=e^{k/\log k}$, $x=x_k=(1+\epsilon)\sigma_{\rho} T_k^{3/2}\sqrt{2\log\log T_k}$. Then for any $\epsilon>0$,
$$ \sum_{k=1}^{\infty}\pr\left(\|\widetilde{\bm C}\|_{T_k}\ge x_k+c_0 T^{3/2}_k\right)\le \sum_{k=1}^{\infty} (k/\log k)^{-(1+\epsilon)}<\infty. $$
Notice $T_{k+1}/T_k\to 1$. It follows that
$$ \limsup_{T\to \infty}\frac{\sup_{0\le t\le T}\|\widetilde{\bm C}_t\|}{T^{3/2}\sqrt{2\log\log T}}
\le \limsup_{k\to\infty} \frac{\|\widetilde{\bm C}\|_{T_k}}{T^{3/2}_k\sqrt{2\log\log T_k}}\le (1+\epsilon)\sigma_{\rho} \; a.s. $$
The "$\le"$ part of \eqref{eq:LILforC} is proven. Next, we consider the "$\ge$" part. It is sufficient to show that, if  $\|\bm b \|=1$,  then
\begin{equation}\label{eq:LILlowbound}
\limsup_{T\to \infty} \frac{ \widetilde{\bm C}_T\bm b^{\prime} }{T^{3/2}\sqrt{2\log\log T}}\ge \sigma_{\rho}\; a.s. 
\end{equation} 
Thus, it is sufficient to consider the one-dimension case. Let $T_k=[e^{k\log k}]$,
$$ D_k=\sigma_1\int_{T_{k-1}}^{T_k}(T_k-s)d   W(s)+\frac{\sigma_2}{1+\rho}\int_{T_{k-1}}^{T_k} \big(T_k^{1+\rho}s^{-\rho}-s\big) d B(t). $$
Then $\{D_k;k\ge 1\}$ is a seqeunce of independent random variables with
\begin{align*}
&C_{T_k}-D_k=\sigma_1\int_0^{T_{k-1}}(T_k-s)d W(t)+\frac{\sigma_2}{1+\rho}\int_{0}^{T_{k-1}} \big(T_k^{1+\rho}s^{-\rho}-s\big) B(t)\\
=& T_kO\big(\sqrt{T_{k-1}\log\log T_{k-1}}\big)+O\big(\sqrt{T_{k-1}^3\log\log T_{k-1}}\big)  
  +T_k^{1+\rho}O\big(\sqrt{T_{k-1}^{1-2\rho}\log\log T_{k-1}}\big)\\
  =&o(\sqrt{T_k^3\log\log T_k})\;\; a.s.  \big(\text{ by } \eqref{eq:proof-CMERW2}),
\end{align*}
and
$$ \Var(D_k)=\sigma_1^2\int_{T_{k-1}}^{T_k}(T_k-s)^2ds+\frac{\sigma_2^2}{(1+\rho)^2}\int_{T_{k-1}}^{T_k} \big(T_k^{1+\rho}s^{-\rho}-s\big)^2ds\sim \sigma^2_{\rho} T_k^3.$$
Then for any $0<\epsilon<1$, when $k$ is large enough,
\begin{align*}
& \pr\left(D_k\ge (1-\epsilon)T_k^{3/2}\sqrt{2 \log\log T_k}\sigma_{\rho}\right)\ge c\frac{1}{x_k/\sqrt{\Var(D_k}}\exp\Big\{-\frac{x_k^2}{2\Var(D_k)}\Big\}\\
\ge & c \exp\{-(1-\epsilon)\log\log T_k\}
=c  (k\log k)^{-(1-\epsilon)}, 
\end{align*}
where $x_k=(1-\epsilon)\sigma_{\rho} T_k^{3/2}\sqrt{2 \log\log T_k}$. By the Borel-Cantelli lemma,
$$\limsup_{k\to \infty}\frac{D_k}{T_k^{3/2}\sqrt{2 \log\log T_k}}\ge (1-\epsilon)\sigma_{\rho}. $$
It follows that
$$\limsup_{k\to \infty}\frac{C_{T_k}}{T_k^{3/2}\sqrt{2 \log\log T_k}}=\limsup_{k\to \infty}\frac{D_k}{T_k^{3/2}\sqrt{2 \log\log T_k}}\ge (1-\epsilon)\sigma_{\rho}. $$
\eqref{eq:LILlowbound} holds.

For (iii),  we let $a_n=\sum_{k=1}^n k^{\rho}$. Then $a_n= \frac{n^{1+\rho}}{1+\rho}\left(1+O(n^{-1})\right)$,
$$ \frac{1}{a_n}\left\|\sum_{k=1}^n\bm T_k-a_n\mu_Z\bm\xi\right\|
\le  \frac{1}{a_n}\sum_{k=1}^nk^{\rho}\left\|\frac{\bm T_k}{k^{\rho}}-\mu_Z\bm\xi\right\|. $$
The results following from Theorem \ref{th:multi-RPWInva2} immediately.

Finally, we consider (iv). In this case, $1/2<\rho<1$. By
\eqref{eq:multi-thRPWInva14}, 
\begin{align*}
 \bm C_n=&\frac{1}{n}\sum_{k=1}^n \left[k^{\rho}\mu_Z\bm \xi+\frac{\sigma_Z}{\sqrt{d}}\bm W(k)+\frac{\mu_Z}{\sqrt{d }}  \bm G_k\right]+o(n^{1/2-\tau}) \; a.s.\\
 =&\frac{\mu_Z}{1+\rho}n^{\rho}\bm \xi+\frac{1}{n}  \left[\frac{\sigma_Z}{\sqrt{d}}\int_0^n \bm W(s)ds+\frac{\mu_Z}{\sqrt{d }}\int_0^n  \bm G_sds\right]+o(n^{1/2-\tau}) \; a.s. 
\end{align*}
\eqref{eq:multi-thCMERWInv3} is proven. Notice \eqref{eq:processidentical}. We still have \eqref{eq:varinceGauss1}-\eqref{eq:varinceGauss2}, and \eqref{eq:LILforC} with $\rho^{\prime}=1-\rho$ taking the place of $\rho$. Hence, by \eqref{eq:multi-thRPWInva14} and  \eqref{eq:multi-thCMERWInv3}, \eqref{eq:multi-thCMERWCLT3} and \eqref{eq:multi-thCMERWLIL3} hold. The proof is now completed.    $\Box$

\subsection{Proof of the Chung LIL}

  We needs some lemmas for the small ball probabilities. 
\begin{lemma}\label{lem:ChungLIL1} (Li \cite[Theorem 1.2]{LiWB1999}) Let $X$ and $Y$ be any two joint Gaussian random vectors in a separable Banach
space with norm $\|\cdot\|$. If
$$\lim_{\epsilon\to 0}(\text{resp.} \; \liminf_{\epsilon\to 0},\; \limsup_{\epsilon\to 0})\epsilon^{\alpha} \log \pr(\|X\|<\epsilon)=-C_X  $$
and
 $$\lim_{\epsilon\to 0}\epsilon^{\alpha} \log \pr(\|Y\|<\epsilon)=0 $$
with $0<\alpha<\infty$ and $0<C_X<\infty$. Then
$$\lim_{\epsilon\to 0}(\text{resp.} \; \liminf_{\epsilon\to 0},\; \limsup_{\epsilon\to 0})\epsilon^{\alpha} \log \pr(\|X+Y\|<\epsilon)=-C_X. $$
\end{lemma}
  
 \begin{lemma}\label{lem:ChungLIL2} (Li \cite[Theorem 1.1]{LiWB1999})  Let $X$   be any centered  Gaussian random vectors in a separable Banach
space $E$.  Then
for any $0 <\lambda < 1$, any symmetric, convex sets $A$ and $B$ in $E$,
$$ \pr( X\in A\cap B)\ge \pr (X\in \lambda A) \pr (X\in (1-\lambda^2)^{1/2} B). $$
\end{lemma}
 
\begin{lemma}\label{lem:ChungLIL3} (Li \cite[Lemma 2.2]{LiWB2001}) If $X_i$, $i=1,\ldots, d$, are independent nonnegative random variables such that
$$\lim_{\epsilon\to 0}\epsilon^{\alpha} \log \pr(X_i<\epsilon)=-C_i,  \;\; i=1,\ldots, d $$
for $0<\alpha<\infty$, then
 $$\lim_{\epsilon\to 0}\epsilon^{\alpha} \log \pr(X_1+\ldots+X_d<\epsilon)=-\left(\sum_{i=1}^d C_i^{1/(1+\alpha)}\right)^{1+\alpha}. $$
\end{lemma}  
  
\begin{lemma}\label{lem:ChungLIL5}  Let $\bm X(t)=(X_1(t),\ldots,X_d(t))$ be a self-similar Gaussian process of index $\tau$, i.e.,
$$ \{\bm X(ct); t\ge 0\}\overset{\mathscr{D}}=\{c^{\tau}\bm X(t); t\ge 0\}, \;\; c>0.  $$
Suppose 
\begin{equation}\label{eq:lemChungLIL5.1}-\log \pr\Big(\sup_{0\le t\le 1} \left\|\bm X(t)\right\|<\epsilon\Big)\preccurlyeq \epsilon^{-\beta}, \;\; \epsilon\to 0  
\end{equation} 
with $0<\beta<\infty$. Then
\begin{equation}\label{eq:lemChungLIL5.2}-\log \pr\left(\sup_{0\le t\le 1} \left\|\int_0^t s^{-\alpha}\bm X(s)ds\right\|<\epsilon\right)\preccurlyeq \epsilon^{-\beta/(\beta+1)}, \;\; \epsilon\to 0, 
\end{equation}
for all $\alpha<\tau+1$.  Furthermore, 
\begin{equation}\label{eq:lemChungLIL5.3}-\log \pr\left(\sup_{0\le t\le 1} \left\|\int_0^t (t^{-\alpha}-s^{-\alpha})\bm X(s)ds\right\|<\epsilon\right)\preccurlyeq \epsilon^{-\beta/(2\beta+1)}, \;\; \epsilon\to 0. 
\end{equation}
and
\begin{equation}\label{eq:lemChungLIL5.4}-\log \pr\left(\sup_{0\le t\le 1} \left\|\int_0^t \bm X(s)ds/t^{\alpha}\right\|<\epsilon\right)\preccurlyeq \epsilon^{-\beta/(\beta+1)}, \;\; \epsilon\to 0  
\end{equation}
for all $\alpha<\tau+1$. 
  Here and in the sequel,  $f(\epsilon)\preccurlyeq g(\epsilon)$ as $\epsilon\to 0$ if $\limsup_{\epsilon\to 0}f(\epsilon)/g(\epsilon)$ is bounded.
\end{lemma} 
{\bf Proof.} For a Gaussian process $\bm Y(t)=(Y_1(t),\ldots,Y_d(t))$, by Lemma \ref{lem:ChungLIL2} we have
\begin{align*}
&\pr\left(\sup_{0\le t\le 1}\|\bm Y(t)\|<\epsilon\right)=\pr\left(\sup_{0\le t\le 1}\sum_{i=1}^d|Y_i(t)|^2<\epsilon^2\right)\\
\ge &\pr\left(\sup_{0\le t\le 1}\sum_{i=1}^{d-1}|Y_i(t)|^2<\epsilon^2/2,\sup_{0\le t\le 1}|Y_d(t)|\le \epsilon/\sqrt{2}\right)\\
\ge &\pr\left(\sup_{0\le t\le 1}\sum_{i=1}^{d-1}|Y_i(t)|^2<\lambda\epsilon^2/2\right)\pr\left(\sup_{0\le t\le 1}|Y_d(t)|\le \sqrt{1-\lambda^2}\epsilon/\sqrt{2}\right). 
\end{align*}
By the induction, it is sufficient to consider the one-dimensional case.  
Let $Y(t)=\int_0^t s^{-\alpha} X(s)ds$. Then $Y(t)$ is a self-similar Gaussian process of index $\tau+1-\alpha$. 
For $0<\delta<1$, let 
$ K(t,s)=s^{-\alpha}I\{\delta\le s\le t\}. $
Then $K(t,s)$ satisfies the H\"older condition
$$ \int_0^1\left|K(t^{\prime},s)- K(t^{\prime\prime},s)\right|ds\le c|t^{\prime}-t^{\prime\prime}|, \;\; t^{\prime}, t^{\prime\prime}\in [0,1], $$
and $\sup_{\delta\le t\le 1}|Y(t)-Y(\delta)|=\sup_{0\le t\le 1}\big|\int_0^1 K(t,s) X(s)ds\big|$.
By \eqref{eq:lemChungLIL5.1} and Theorem 6.1 of Li and Linde \cite{LiLinde1998}, we have
$$ -\log \pr\Big(\sup_{\delta\le t\le 1}|Y(t)-Y(\delta)|< \epsilon\Big)\preccurlyeq \epsilon^{-\frac{\beta}{\beta+1}}. $$
On the other hand, it is easily checked that 
$$ \lim_{\epsilon\to 0}\epsilon^{\gamma}\log \pr\Big(|Y(\delta)|< \sqrt{1-\lambda^2}\epsilon\Big)=0, \forall \gamma>0. $$
By Lemma \ref{lem:ChungLIL2}, it follows that 
$$-\log \pr\Big(\sup_{\delta\le t\le 1}|Y(t)|< \epsilon\Big)
\preccurlyeq \epsilon^{-\frac{\beta}{\beta+1}}. 
$$
 Thus, there exist two  positive constants  $\kappa$ and $\epsilon_0$ such that
 \begin{align*}
  \log \pr\Big(\sup_{\delta\le t\le 1} |Y(t)|< \epsilon\Big) 
\ge
 -\kappa\epsilon^{-\frac{\beta}{\beta+1}}, \;\; 0<\epsilon\le \epsilon_0.
 \end{align*}
  Choose $0<\delta<1$ such that  $\lambda=2\delta^{1+\tau-\alpha}<1$. By  Lemma \ref{lem:ChungLIL2} again, 
 \begin{align*}
 & \pr\Big(\sup_{0\le t \le 1} \left|Y(t)\right|<\epsilon\Big)
\ge   \pr\Big(\sup_{0\le t\le \delta} \left|Y(t)\right|< \lambda\epsilon\Big)\pr\Big( \sup_{\delta \le t\le 1} \left|Y(t)\right|< \sqrt{1-\lambda^2}\epsilon \Big)\\
= &   \pr\Big(\sup_{0\le t\le 1}\left|Y(t)\right|<\lambda\delta^{\alpha-1-\tau} \epsilon\Big)\pr\Big( \sup_{\delta\le t\le 1} \left|Y(t)\right|< \sqrt{1-\lambda^2}\epsilon \Big).
 \end{align*}
 For $0<\epsilon\le \epsilon_0/8$, choose $k$ such that $2^k\epsilon\le \epsilon_0<2^{k+1}\epsilon$. Then, we have 
 \begin{align*}
&\log\pr\Big(\sup_{0\le t \le 1} \left|Y(t)\right|<\epsilon\Big)    
\ge  \log\pr\Big(\sup_{0\le t \le 1} \left|Y(t)\right|<2\epsilon\Big)-\kappa\big(\sqrt{1-\lambda^2} \epsilon\big)^{-\beta/(\beta+1)}\\
\ge &  \log\pr\Big(\sup_{0\le t \le 1} \left|Y(t)\right|<4\epsilon\Big)-\kappa\big(\sqrt{1-\lambda^2} \epsilon\big)^{-\beta/(\beta+1)}
-\kappa(\sqrt{1-\lambda^2} 2\epsilon)^{-\beta/(\beta+1)}\\
\ge & \cdots\ge   \log \pr\Big(\sup_{0\le t \le 1} \left|Y(t)\right|<2^{k+1}\epsilon\Big)-\epsilon^{-\frac{\beta}{\beta+1}}\frac{\kappa}{(1-\lambda^2)^{\beta/(2\beta+2)}} 
\sum_{i=0}^k 2^{-k\beta/(\beta+1)}\\
\ge  &  \log\pr\Big(\sup_{0\le t \le 1} \left|Y(t)\right|<\epsilon_0\Big)-c\epsilon^{-\frac{\beta}{\beta+1}}.
\end{align*}  
\eqref{eq:lemChungLIL5.2}  is proven.

Write $X^{\ast}(t)=\int_0^t X(s)ds$. Then $X^{\ast}(t)$ is a self-similar Gaussian process of index $\tau+1$,  
$$ \int_0^t(t^{-\alpha}-s^{-\alpha})X(s)ds=-\alpha\int_0^t s^{-\alpha-1}X^{\ast}(s) ds, $$
and 
$$-\log \pr\Big(\sup_{0\le t\le 1}|X^{\ast}(t)|<\epsilon\Big)\preccurlyeq \epsilon^{-\beta/(\beta+1)}, \; \epsilon\to 0$$
by \eqref{eq:lemChungLIL5.2}. With $X^{\ast}(t)$, $\beta/(\beta+1)$, and $\alpha+1$ taking the places of $X(t)$, $\beta$, and $\alpha$ respectively, we obtain  \eqref{eq:lemChungLIL5.3}. \eqref{eq:lemChungLIL5.4} follows from \eqref{eq:lemChungLIL5.2} and \eqref{eq:lemChungLIL5.3} by Lemma \ref{lem:ChungLIL1}, and the proof of the lemma is completed. $\Box$

We first consider the Chung type law of the iterated logarithm of $\bm T_n$.

{\bf Proof of Proposition \ref{prop:smallball1}.} When $\rho_1=\rho_2=0$, \eqref{eq:propsmallball1.2} is obvious since
\begin{equation}\label{eq:proof-propsmallball1.0} \{\bm I(t),t\ge 0\}\overset{\mathscr{D}}=\{(\sigma_1^2+\sigma_2^2)^{1/2} \bm B(t);t\ge 1\}. 
\end{equation}
Hence, by Lemma \ref{lem:ChungLIL1},  for \eqref{eq:propsmallball1.2} it is sufficient to show that 
\begin{equation}\label{eq:proof-propsmallball1.1} -\log \pr\left(\sup_{0\le t\le 1}\left\|\int_0^t\big\{(t/s)^{\gamma}-1\}d\bm B(s)\right\|<\epsilon\right)\preccurlyeq \epsilon^{-2/3}, \; \epsilon\to 0
\end{equation}
for $\gamma<1/2$ and $\gamma\ne 0$.  Let 
\begin{equation}\label{eq:proof-propsmallball1.2} \bm Y(t)=\int_0^ts^{-\gamma} d\bm B(s)\; \text{ and } \bm X(t)=\sqrt{1-2\gamma} \bm Y(t^{1/(1-2\gamma)}). 
\end{equation}
 Then, $\bm X(t)$ is a $d$-dimensional  standard Browian motion, and $d \bm B(s)=s^{\gamma}d\bm Y(s)$. Thus,
\begin{align*}
&\sup_{0\le t\le 1}\left\|\int_0^t\big\{(t/s)^{\gamma}-1\}d \bm B(s)\right\|
=\sup_{0\le t\le 1}\left\|\int_0^t\big\{t^{\gamma}-s^{\gamma}\}d \bm Y(s)\right\|
=\sup_{0\le t\le 1}\left\|\gamma\int_0^ts^{\gamma-1}\bm Y(s)ds\right\|\\
& \;\; = \sup_{0\le t\le 1}\left\|\frac{\gamma}{(1-2\gamma)^{3/2}}\int_0^{t^{1-2\gamma}}s^{\frac{3\gamma-1}{1-2\gamma}}\bm X(s)ds\right\|
=\frac{|\gamma|}{(1-2\gamma)^{3/2}}\sup_{0\le t\le 1}\left\|\int_0^ts^{\frac{3\gamma-1}{1-2\gamma}}\bm X(s)ds\right\|.
\end{align*}
By Lemma \ref{lem:ChungLIL5} with $\tau=1/2$, $\beta=2$, and $\alpha=\frac{1-3\gamma}{1-2\gamma}$, \eqref{eq:proof-propsmallball1.1} holds since $\frac{1-3\gamma}{1-2\gamma}<3/2$.  The proof of Proposition \ref{prop:smallball1} is completed. $\Box$

{\bf Proof of Corollary \ref{cor:CLIL}.} We consider the Gaussian process $\widetilde{\bm T}_t=(\sigma_Z\bm W(t)+\mu_Z\bm G_t)/\sqrt{d}$. Suppose $\rho<1/2$. By Proposition \ref{prop:smallball1}, 
$$ \lim_{\epsilon\to 0}\epsilon^2\log \pr\left(\sup_{0\le t\le 1} \|\widetilde{\bm T}_t\|<\epsilon\right)=-\frac{j_{(d-2)/2}^2}{2d}\ep[Z_1^2]. $$
Notice that $\{\bm T_{Tt}/\sqrt{T};t\ge 1\}\overset{\mathscr{D}}=\{\bm T_t;t\ge 0\}$. 
It follows from  standard
arguments along with an application of the Borel-Cantelli lemma (see for example Cs\"org\H o and R\'ev\'esz \cite[Chap. 1]{CR1981}) that
$$\liminf_{T\to \infty}\Big(\frac{\log\log T}{T}\Big)^{1/2} \max_{0\le t\le T}\|\widetilde{\bm T}_t\|=\frac{j_{(d-2)/2}}{\sqrt{2d}}(\ep[Z_1^2])^{1/2}\; a.s.  $$ 
  The proof of \eqref{eq:CLIL2} is completed by \eqref{eq:multi-thRPWInva2}. 
 
 When $1/2<\rho<1$, notice that
 \begin{equation}\label{eq:processidentical}\left\{\bm G_t;t\ge 0\right\}\overset{\mathscr{D}}=\left\{\frac{t^{1-\rho}}{\sqrt{2\rho-1}}\bm B(t^{2\rho-1});t\ge 0\right\}\overset{\mathscr{D}}=\left\{\int_0^t(t/s)^{1-\rho}d\bm B(s);t\ge 0\right\}. 
 \end{equation}
 By  \eqref{eq:multi-thRPWInva14}, the proof is the same as that in the case of $\rho<1/2$. 
 $\Box$
 
Next, we consider the Chung type law of the iterated logarithm of the center of mass $\bm C_n$. 

{\bf Proof of Proposition \ref{prop:smallball2}.} First, the existence of the limit $\kappa_d$ in \eqref{eq:SBPofIB} can be found in Khoshnevisan and  Shi \cite{KhoshnevisanShi1998} based on the subadditivity. When $d=1$, Chen and Li \cite{ChenLi2003} have shown that $3/8\le \kappa_d\le (2\pi)^{2/3}\cdot(3/8)$.  Notice that
\begin{align*}
&\pr\left(\sup_{0\le t\le 1}\Big|\int_0^t B_1(s)ds\Big|< \epsilon\right)  
\ge \pr\left(\sup_{0\le t\le 1}\Big\|\int_0^t\bm B(s)ds\Big\|< \epsilon\right)   \\
& \qquad \ge \pr\left(\sum_{i=1}^d \sup_{0\le t\le 1}\Big|\int_0^tB_i(s)ds\Big|^2< \epsilon^2\right).
\end{align*}  
By Lemma \ref{lem:ChungLIL3}, $\kappa_1\le \kappa_d\le d^{4/3}\kappa_1$. 

Similar to Proposition \ref{prop:smallball1}, by Lemma \ref{lem:ChungLIL1} and noticing \eqref{eq:proof-propsmallball1.0}, for the second equality in \eqref{eq:propsmallball2.2} it is sufficient to show that for $\alpha<3/2$, 
\begin{equation}\label{eq:weightedSBIN}
\lim_{\epsilon\to 0}\epsilon^{2/3}\log \pr\left(\sup_{0< t\le 1}\frac{1}{t^{\alpha}}\Big\|\int_0^t\bm B(s)ds\Big\|< \epsilon\right)  
=-\frac{\kappa_d}{1-2\alpha/3}
\end{equation}
and
\begin{equation}\label{eq:proofprop2.22}  -\log \pr\left(\sup_{0<t\le 1}\big|\bm J_{\gamma}(t)/t^{\alpha}\big|< \epsilon\right)\preccurlyeq \epsilon^{-2/5},\; \epsilon\to 0,  
\end{equation}
where $\bm G_{\gamma}(t)=\int_0^t \{(t/s)^{\gamma}-1\} d \bm B(s)$ and $\bm J_{\gamma}(t)=:\int_0^t \bm G_{\gamma}(s) ds$, $\gamma<1/2$.

Notice that $\bm G_{\gamma}(t)$ is a self-similar Gaussian process with index $1/2$.  
 \eqref{eq:proofprop2.22} follows from \eqref{eq:proof-propsmallball1.1} by Lemma \ref{lem:ChungLIL5}.
 
 For \eqref{eq:weightedSBIN}, we 
write $P(\epsilon)=\log \pr\left(\sup_{0\le t\le 1}\frac{1}{t^{\alpha}}\Big\|\int_0^t\bm B(s)ds\Big\|< \epsilon\right)$. Notice $$\lim_{\epsilon\to 0}\epsilon^2\log \pr\Big(\sup_{0\le t\le 1}\big\|\bm B(t)\big\|< \epsilon\Big)=-\frac{j_{(d-2)/2}^2}{2}.$$ By Lemma \ref{lem:ChungLIL5},  we have
$0\le -P(\epsilon)\preccurlyeq \epsilon^{-\frac{2}{2+1}}=\epsilon^{-2/3}, \epsilon\to 0. $
Thus, $\limsup_{\epsilon\to 0}\epsilon^{2/3}P(\epsilon)$ and $\liminf_{\epsilon\to 0}\epsilon^{2/3}P(\epsilon)$ are finite. 
Now, for $0<\delta<1$ we have
\begin{align*}
&\pr\left(\sup_{0< t\le 1}\frac{1}{t^{\alpha}}\Big\|\int_0^t\bm B(s)ds\Big\|< \epsilon\right)\\
\ge &\pr\left(\sup_{0< t\le \delta}\frac{1}{t^{\alpha}}\Big\|\int_0^t\bm B(s)ds \Big\|< \epsilon, 
\sup_{\delta\le t\le 1}\Big\|\int_0^t\bm B(s)ds\Big\|<  \delta^{|\alpha|}\epsilon\right)\\  
= &\pr\left(\sup_{0< t\le \delta}\frac{1}{t^{\alpha}}\Big\|\int_0^t\bm B(s)ds\Big\|< \epsilon, \right. \\
&\qquad \left.\sup_{\delta\le t\le 1}\Big\|\int_{\delta}^t(\bm B(s)-\bm B(\delta))ds+ (t-\delta)\bm B(\delta)+\int_0^{\delta} \bm B(s) ds\Big\|<\delta^{|\alpha|}\epsilon\right).
\end{align*}
It is easily checked that
$$ \lim_{\epsilon\to 0}\epsilon^{2/3}\log\pr\left(\sup_{\delta\le t\le 1}\Big\| (t-\delta)\bm B(\delta)+\int_0^{\delta} \bm B(s) ds\Big\|<\epsilon\right)=0. $$
By Lemma \ref{lem:ChungLIL1} and the independence,
\begin{align*}
&\liminf_{\epsilon\to 0}\epsilon^{2/3}P(\epsilon)\\
\ge &\liminf_{\epsilon\to 0}\epsilon^{2/3}\log \pr\left(\sup_{0< t\le \delta}\frac{1}{t^{\alpha}}\Big\|\int_0^t\bm B(s)ds\Big\|< \epsilon, \sup_{\delta\le t\le 1}\Big\|\int_{\delta}^t(\bm B(s)-\bm B(\delta))ds \Big\|<\delta^{|\alpha|} \epsilon\right)\\
= &\liminf_{\epsilon\to 0}\epsilon^{2/3}\log\left\{\pr\left(\sup_{0< t\le \delta}\frac{1}{t^{\alpha}}\Big\|\int_0^t\bm B(s)ds\Big\|< \epsilon\right) \pr\left( \sup_{\delta\le t\le 1}\Big\|\int_{\delta}^t(\bm B(s)-\bm B(\delta))ds \Big\|< \delta^{|\alpha|}\epsilon\right)\right\}\\
=& \liminf_{\epsilon\to 0}\epsilon^{2/3}P(\delta^{\alpha-3/2}\epsilon)
+\lim_{\epsilon\to 0}\epsilon^{2/3} \log\pr\left( \sup_{0\le t\le 1}\Big\|\int_{0}^t \bm B(s)ds \Big\|< \frac{\delta^{|\alpha|}\epsilon}{(1-\delta)^{3/2}}\right)  \\
=& \delta^{1-2\alpha/3}\liminf_{\epsilon\to 0}\epsilon^{2/3}P(\epsilon)
-\kappa_d\frac{1-\delta}{\delta^{2|\alpha|/3}}.
\end{align*}
 We conclude that
$$ \liminf_{\epsilon\to 0}\epsilon^{2/3}P(\epsilon)\ge -\kappa_d\frac{1-\delta}{\delta^{2|\alpha|/3}(1-\delta^{1-2\alpha/3})}\to -\frac{\kappa_d}{1-2\alpha/3}
\text{
as }  \delta\to 1. $$ 

On the other hand, 
for $0<\delta<1$,  
\begin{align*}
&\pr\left(\sup_{0< t\le 1}\frac{1}{t^{\alpha}}\Big\|\int_0^t\bm B(s)ds\Big\|< \epsilon\right) \\
\le &\pr\left(\sup_{0< t\le \delta}\frac{1}{t^{\alpha}}\Big\|\int_0^t\bm B(s)ds\Big\|< \epsilon, 
\sup_{\delta\le t\le 1}\Big\|\int_0^t\bm B(s)ds\Big\|< \delta^{-|\alpha|} \epsilon\right)\\
= &\pr\left(\sup_{0< t\le \delta}\frac{1}{t^{\alpha}}\Big\|\int_0^t\bm B(s)ds\Big\|< \epsilon, \right. \\
&\qquad \left.\sup_{\delta\le t\le 1}\Big\|\int_{\delta}^t(\bm B(s)-\bm B(\delta))ds+ (t-\delta)\bm B(\delta)+\int_0^{\delta} \bm B(s) ds\Big\|< \delta^{-|\alpha|}\epsilon\right).
\end{align*}
With the same argument above, we have
\begin{align*}
& \limsup_{\epsilon\to 0}\epsilon^{2/3}P(\epsilon)\\
\le & \limsup_{\epsilon\to 0}\epsilon^{2/3}P(\delta^{\alpha-3/2}\epsilon)+
\lim_{\epsilon\to 0}\epsilon^{2/3} \log\pr\left( \sup_{0\le t\le 1}\Big\|\int_{0}^t \bm B(s)ds \Big\|< \frac{\delta^{-|\alpha|}\epsilon}{(1-\delta)^{3/2}}\right)\\
=& \delta^{1-2\alpha/3}\limsup_{\epsilon\to 0}\epsilon^{2/3}P(\epsilon)- \kappa_d(1-\delta)\delta^{2|\alpha|/3}.
\end{align*} 
Thus
$$ \limsup_{\epsilon\to 0}\epsilon^{2/3}P(\epsilon)\le -\kappa_d
\frac{(1-\delta)\delta^{2|\alpha|/3}}{1-\delta^{1-2\alpha/3}} \to -\frac{\kappa_d}{1-2\alpha/3}
\;\text{
as } \; \delta\to 1. $$ 
\eqref{eq:weightedSBIN} is proven, and the proof of the second equality in \eqref{eq:propsmallball2.2}    is completed. 

Notice that $\bm I(t)$ is a self-similar Gaussian process of index $1/2$ and \eqref{eq:propsmallball1.2} holds. By Lemma \ref{lem:ChungLIL5},
$$-\log  \pr\left(\sup_{0< t\le 1}\Big\| \int_0^t  \big(t^{-\alpha}-s^{-\alpha}\big)\bm I(s)ds\Big\|<\epsilon\right)\preccurlyeq
\epsilon^{-\frac{2}{2\ast 2+1}}=\epsilon^{-2/5}, \;\; \alpha<3/2. $$ 
By Lemma \ref{lem:ChungLIL1}, the first equality in \eqref{eq:propsmallball2.2} holds.

Finally, \eqref{eq:propsmallball2.3} follows from the estimates given in \eqref{eq:propsmallball2.2} and standard
arguments along with an application of the Borel-Cantelli lemma. $\Box$ 
 
 {\bf Proof of Corollary \ref{cor:CLIL2}.}
 Let $\widetilde{\bm C}_T$ be defined as in \eqref{eq:CMEWGauss1} with $\sigma_1=\sigma_Z/\sqrt{d}$ and $\sigma_2=\mu_Z/\sqrt{d}$. By noticing \eqref{eq:processidentical}, it is sufficient to consider the case of  $\rho<1/2$.  By Proposition \ref{prop:smallball2} with $\rho_1=0$, $\rho_2=\rho$ and $\alpha=1$, 
$$\liminf_{T\to \infty}\frac{(\log\log T)^{3/2}}{T^{1/2}}  \max_{0< t\le T}\|\widetilde{\bm C}_t/t\|=(3\kappa_d)^{3/2} (\ep Z_1^2/d)^{1/2}\; a.s.  $$ 
The proof of \eqref{eq:CLIL2.1} is completed by \eqref{eq:multi-thCMERWInv1}. 
 $\Box$

\subsection{Proof of the almost sure CLT}

{\bf Proof of Corollary \ref{cor:ASCLT}.} We first show (i).  Suppose $\rho<1/2$ at first. 
 Let $\widetilde{\bm T}_T$ and $\widetilde{\bm C}_T$ be defined as in \eqref{eq:CMEWGauss1},  and $\bm X(t)=(  \widetilde{\bm T}_{t}/\sqrt{t},   \widetilde{\bm C}_{t}/ t^{3/2})$, $t>0$. 
It is easily check that
$$\Cov\left\{\frac{\widetilde{\bm T}_S}{\sqrt{S}}, \frac{\widetilde{\bm T}_T}{\sqrt{T}}\right\}
=\left[\sigma_1^2\big(S/T)^{1/2}+\frac{\sigma_2^2}{1-2\rho} \big(S/T\big)^{1/2-\rho}\right] I_d,\; S\le T. $$
Let  $\bm U(t)=\bm X(e^t)=\big(\bm U_1(t), \bm U_2(t)\big)$, $t\in \mathbb R$,  where  $\bm U_1(t)=\frac{\widetilde{\bm T}_{e^t}}{\sqrt{e^t}}$, $\bm U_2(t)=\frac{\widetilde{\bm C}_{e^t}}{e^{3t/2}}$.  Then $\bm U(0)=\bm X(1)=(\widetilde{\bm T}_1,\widetilde{\bm C}_1)$, $\bm U_1(t)$ is a mean zero Gaussian process with
$$ \Cov\{\bm U_1(t),\bm U_1(s)\}=\left[\sigma_1^2\exp\{-\frac{1}{2}|t-s|\}+\frac{\sigma_2^2}{1-2\rho}\exp\{-\big(\frac{1}{2}-\rho\big)|t-s|\}\right]I_d. $$
It follows that $\bm U_1(t)$ is a mean zero stationary $d$-dimensional  Gaussian process.
Notice that
$$ \bm U_2(t)=e^{-3t/2} \int_0^{e^t} \widetilde{\bm T}_x dx=\int_{-\infty}^0 \bm U_1(t+y)e^{3y/2}dy. $$ 
Write $\bm O_t(y)=(\bm U_1(t), \bm U_1(t+y)e^{3y/2})$. Then $\bm O_t(\cdot)$ is a function-valued stationary process. 
Hence, $\bm U(t)$ is a $2d$-dimensional stationary   Gaussian process.  By the Birkhoff ergodic theorem \cite{Birkhoff1931},
\begin{align} \label{eq:ASCLTG1} 
\overline{f}= &\lim_{T\to \infty}\frac{1}{\log T}\int_1^{T} \frac{f\big(\bm X(t)\big)}{t}  
=   \lim_{\tau\to \infty}\frac{1}{\tau}\int_0^{\tau} f\big(\bm U(t)\big)\; \text{ exists almost surely}   \\
& \; \text{ whenever }\; \ep\big|f\big(\bm U(0)\big)\big|=\ep\big|f\big(\bm X(1)\big)\big|<\infty \nonumber
\end{align}
and $\ep[\overline{f}]=\ep\big[f\big(\bm X(1)\big)\big]$.

Let $C_1>0$ such that $C_1\psi(A_0)\ge C_0$, define $\psi^{\ast}(x)=C_1\psi(A_0)$, $0\le x\le A_0$,  $\psi^{\ast}(x)=C_1\psi(x)$, $ x>A_0$. Then $|f(\bm x)|\le \psi^{\ast}(\|\bm x\bm M\|)$, $\psi^{\ast}(x)\exp\{-H_0x^2\}$ is non-increasing on $[0,\infty)$. Thus, without loss of generality, we assume $A_0=0$. Let $\mathcal C_b$, $\mathcal C_{b,Lip}$ be the sets of bounded continuous functions, bounded  Lipschitz functions on $\mathbb R^{2d}$, respectively.  We will finish the proof of (i) by four steps. 

{\em Step 1}. Let $f(\bm x)=\big(\psi(\|\bm x\bm M\|)-A)^+$, $\|\bm x\|$  or $f\in \mathcal C_{b,Lip}$. We have that
 \begin{align}  \label{eq:ASCLTG2} \lim_{T\to \infty}& \frac{1}{\log T}   \int_1^{T} \frac{f\big(\bm X(t)+\bm\delta_t\big)}{t}=\overline{f} \; a.s. \\
&\; \text{ whenever }\; \bm\delta_t=o(1/\sqrt{\log\log t})\; a.s. \nonumber \\
  \text{ and }&  \overline{f}=\ep\big[f\big(\bm X(1)\big)\big]\; a.s. \label{eq:ASCLTG3}
\end{align}

If $\bm f\in \mathcal{C}_{b,Lip}$ or $\|\bm x\|$,  then $f(\bm X(t)+\bm\delta_t)-f(\bm X(t))=O(\|\bm\delta_t\|)=o(1)$ a.s. \eqref{eq:ASCLTG2} holds obviously. 

 Let $f(\bm x)=\big(\psi(\|\bm x\bm M\|)-A)^+$.   By the property that $\psi(x)\exp\{-H_0x^2\}$ is non-increasing on $[0,\infty)$, we have $\psi(x+\delta)\le \psi(x)\exp\{H_0( 2x\delta+\delta^2)\}$ for $x\ge 0$ and $\delta\ge 0$. Then
 \begin{align*} \psi\left(\|(\bm X(t)+\bm\delta_t)\bm M\|\right)\le &
 \psi(\|\bm X(t)\bm M\|)\exp\left\{H_0(2\|\bm X(t)\bm M\|\cdot \|\bm\delta_t\|+\|\bm \delta_t\bm M\|^2)\right\} \\
 = &\psi(\|\bm X(t)\bm M\|)\exp\{o(1)\}\;\;a.s.
 \end{align*}
by noticing that $\bm X(t)=O\big((t\log\log t)^{1/2}\big)$ a.s. and $\bm\delta_t=o(1/\sqrt{\log\log t})$ a.s. Similarly 
$ \psi\left(\|(\bm X(t)+\bm\delta_t)\bm M\|\right)\ge 
  \psi(\|\bm X(t)\bm M\|)\exp\{o(1)\}\;\; a.s.$ 
 Thus
 $$\psi\left(\|(\bm X(t)+\bm\delta_t)\bm M\|\right) 
 = \psi(\|\bm X(t)\bm M\|)\exp\{o(1)\}\;\; a.s.$$
It follows that
$$\left(\psi\left(\|(\bm X(t)+\bm\delta_t)\bm M\|\right) -A\right)^+
 = \left(\psi(\|\bm X(t)\bm M\|)-A\right)^+(1+o(1))+o(1)\;\; a.s.$$
 Hence, \eqref{eq:ASCLTG2} holds for $f(\bm x)=\big(\psi(\|\bm x\bm M\|)-A)^+$. 
  
For any given $t_0\ge 0$, define $\bm B^{\ast}(t)=\bm W^{\ast}(t)=\bm 0$ if $t\le t_0$, $\bm B^{\ast}(t)=\bm B(t)-\bm B(t_0)$ and 
$\bm W^{\ast}(t)=\bm W(t)-\bm W(t_0)$ if $t\ge t_0$. With $\bm B^{\ast}(t)$ and $ \bm W^{\ast}(t)$ taking the place of $\bm B(t)$ and $\bm W(t)$ in \eqref{eq:CMEWGauss2}, we redefine $\widetilde{\bm T}_T$, $\widetilde{\bm C}_T$ and $\bm X(t)$, and denote them by 
$\widetilde{\bm T}_T^{\ast}$, $\widetilde{\bm C}_T^{\ast}$ and $\bm X^{\ast}(t) $, respectively. It is easily seen that
$$\bm\delta_t= \bm X^{\ast}(t)-\bm X(t)=O(t^{-1/2}+t^{-(1/2-\rho)})=o(t^{-\tau}) \; a.s. \text{ as } t\to\infty. $$
Thus
$$ \overline{f}=\lim_{T\to \infty}  \frac{1}{\log T}   \int_1^{T} \frac{f\big(\bm X^{\ast}(t)\big)}{t}\; a.s.  $$
by \eqref{eq:ASCLTG2}. It follows that $\overline{f}$ is independent of $\bm B(t_0),\bm W(t_0)$. By the arbitrariness of $t_0\ge 0$, $\overline{f}$ is independent of $\sigma(\bm B(t),\bm W(t);t\ge 0)$. Hence,  $\overline{f}$ is independent of itself and must be a constant. \eqref{eq:ASCLTG3} is proven.

{\em Step 2}.  We show that \eqref{eq:ASCLTG2} and \eqref{eq:ASCLTG3} hold for all almost everywhere continuous functions $f(\bm x)$  which satisfy \eqref{eq:ASCLTcond1} and \eqref{eq:ASCLTcond2}, and, moreover, the exceptional set of probability $0$ can be chosen universally for all such $f$. 

Let
$$\mathcal{Q}_T^{\omega}=\frac{1}{\log T}\int_1^T\frac{1}{t}\Delta_{\bm X(t,\omega)+\bm\delta_t(\omega)}, $$
be the empirical measure associated with the process $\{\bm X(t)+\bm \delta_t\}$, where $\Delta_{\bm x}$ is the probability measure on $\mathbb R^{2d}$ which assigns
its total mass to $\bm x$. Denote $\mathcal{Q}$ to be the measure generated by $\widetilde{\bm C}$, $\mathcal{Q}(A)=\pr(\widetilde{\bm C}\in A)$. By Step 1, 
\begin{equation}\label{eq:empiricalmeasure}\lim_{T\to \infty}\int f d\mathcal{Q}_T^{\omega} =\int f d\mathcal{Q}  
\end{equation} 
almost surely for all $f\in \mathcal C_{b,Lip}$,   $f(\bm x)=\|\bm x\|$  and $f(\bm x)=\big(\psi(\|\bm x \bm M\|)-A\big)^+$, $A\ge 0$. Notice that for a continuous function $f$ with $\|f\|=\sup_{\bm x}|f(\bm x)|\le 1$ and any $\epsilon>0$, $A>\ep[\|\widetilde{\bm C}\|]/\epsilon$, there exists a bounded Lipschitz function $g_{\epsilon}$ with $\|g_{\epsilon}\|\le 1$ such that
$$ \big|f(\bm x)-g_{\epsilon}(\bm x)\big|I\{\|\bm x\|\le A\}\le \epsilon. $$
Thus
$$ \big|f(\bm x)-g_{\epsilon}(\bm x)\big|\le \epsilon +2\frac{\|\bm x\|\wedge A}{A}, $$
which implies
$$ \limsup_{T\to \infty}\int |f-g_{\epsilon}|d\mathcal{Q}_T\le \epsilon+\frac{\ep[\|\widetilde{\bm C}\|\wedge A]}{A} <2\epsilon\;\; a.s.$$ 
and $\ep[|f(\widetilde{\bm C})-g_{\epsilon}(\widetilde{\bm C})|]<2\epsilon$. 
It follows that \eqref{eq:empiricalmeasure} also holds almost surely  if $f\in \mathcal{C}_b$.

Since $\mathcal C_b$ with the norm $\|f\|=\sum_{k=1}^{\infty}\sup\limits_{\|\bm x\|\le k}|f(\bm x)|/2^k$ is a complete and separable space and so there is a countable dense subset, 
it can be shown that there exists a null event $\Omega_0$ such that \eqref{eq:empiricalmeasure} holds for all $f\in \mathcal C_b$, $f(\bm x)=\|\bm x\|$  and $f(\bm x)=\big(\psi(\|\bm x\bm M\|)-A\big)^+$, $A\ge 0$, when $\omega\notin \Omega_0$. 

Now, assume that $f(\bm x)$ is an almost everywhere continuous function which satisfies \eqref{eq:ASCLTcond1} and \eqref{eq:ASCLTcond2}, and, $\omega\notin \Omega_0$.  Write $\psi(\bm x)=\psi(\|\bm x\bm M\|)$ and  $f^{(A)}(\bm x)=(-A)\vee f(\bm x)\wedge A$.
Then $|f-f^{(A)}|=(|f|-A)^+ \le (\psi-A)^+$. It follows that
$$ \limsup_{T\to\infty}\int |f-f^{(A)}|d \mathcal{Q}_T^{\omega}\le \limsup_{T\to\infty}\int  (\psi-A)^+d \mathcal{Q}_T^{\omega}\le \int(\psi-A)^+ d\mathcal{Q}\to 0 \text{ as } A\to \infty. $$

For $f^{(A)}$,  $0\le f^{(A)}+A\le 2A$.   Recall that  $D_f$ is the set of discontinuities of $f$. Let $F_s=\{\bm x: f(\bm x)+A\ge s\}$. Then $F_s\subset cl(F_s)\subset F_s\cup D_f$, $\mathcal{Q}(D_f)=\pr(\widetilde{\bm C}\in D_f)=0$, where $cl(F_s)$ is the closure of $F_s$.   Thus, $\mathcal{Q}(F_s)=\mathcal{Q}(cl(F_s))$.  
  We have
\begin{align*}
   \int  (f^{(A)}+A) d\mathcal{Q}_T^{\omega}=&   \int_0^{\infty} \mathcal{Q}_T^{\omega}(F_s)ds  
\le     \int_0^{3A} \mathcal{Q}_T^{\omega}(F_s)ds  
\le   \int_0^{3A} \mathcal{Q}_T^{\omega}(cl(F_s))ds.  
\end{align*}
By noting   that $d(\bm x, cl(F_s))\wedge \epsilon$ is a bounded continuous function and $d(\bm x, cl(F_s))=0\Longleftrightarrow \bm x \in cl(F_s)$, we have
\begin{align*}
&\limsup_{T\to \infty}  \int  (f^{(A)}+A) d\mathcal{Q}_T^{\omega} 
\le     \limsup_{T\to \infty} \int_0^{3A}\left[\int \Big(1-\frac{d(\cdot,cl(F_s))\wedge \epsilon}{\epsilon}\Big) d \mathcal{Q}_T^{\omega}\right]ds \\
=&\int_0^{3A}\left[\lim_{T\to \infty}\int \Big(1-\frac{d(\cdot,cl(F_s))\wedge \epsilon}{\epsilon}\Big) d \mathcal{Q}_T^{\omega}\right]ds 
= \int_0^{3A}\left[ \int \Big(1-\frac{d(\cdot,cl(F_s))\wedge \epsilon}{\epsilon}\Big) d \mathcal{Q}\right]ds\\
&\overset{\epsilon\to 0}\longrightarrow   \int_0^{3A} \mathcal{Q}(cl(F_s)) ds=\int_0^{3A} \mathcal{Q}(F_s) ds\le \int  (f^{(A)}+A) d\mathcal{Q}.
\end{align*}
It follows that
$$ \limsup_{T\to \infty} \int  f^{(A)}d\mathcal{Q}_T^{\omega}\le \int f^{(A)} d\mathcal{Q}\to \int f  d\mathcal{Q}, $$
as $A\to \infty$. Hence,
$$ \limsup_{T\to \infty} \int  f d\mathcal{Q}_T^{\omega}\le   \int f  d\mathcal{Q}. $$ 
For $-f$, we have a similar inequality. Thus \eqref{eq:empiricalmeasure} holds. 

{\em Step 3}. We show \eqref{eq:ASCLT1.1}.

Let $\sigma_1=\sigma_Z/\sqrt{d}$, $\sigma_2=\mu_Z/\sqrt{d}$. Then 
$$ \bm T_k/\sqrt{k}-\widetilde{\bm T}_k/\sqrt{k}=o(k^{-\tau}) \text{ and } \bm C_k/\sqrt{k}-\widetilde{\bm C}_k/k^{3/2}=o(k^{-\tau})\; a.s. $$
by \eqref{eq:multi-thRPWInva2} and \eqref{eq:multi-thCMERWInv1}.  Write 
$$\bm X^{\ast}(t)= (\bm T_k/\sqrt{k}, \bm C_k/\sqrt{k}), \; t\in (k-1,k]. $$
It can be check that
\begin{align*} &\sup_{k-1\le t\le k}\left\|\widetilde{\bm T}_t/\sqrt{t}-\widetilde{\bm T}_k/\sqrt{k}\right\|=o(k^{-1/3})\;\; a.s.  
 \text{ and } \\
 & \sup_{k-1\le t\le k}\left\|\widetilde{\bm C}_t/t^{3/2}-\widetilde{\bm C}_k/k^{3/2}\right\|=o(k^{-1/3})\;\; a.s. 
\end{align*}
It follows that
$$ \bm X^{\ast}(t)-\bm X(t)=o(t^{-\tau})\; a.s. \; \text{ as } t\to \infty. $$
Hence, by \eqref{eq:ASCLTG2},
\begin{align*}
 &\frac{1}{\log n}\sum_{k=2}^n \frac{f\big((\bm T_k,\bm C_k)/\sqrt{k}\big)}{k} 
 = \frac{1}{\log n}\sum_{k=2}^n\int_{k-1}^k \frac{f\big( \bm X^{\ast}(t)\big)}{k}\\
 =& \frac{1}{\log n}\int_{1}^n \frac{f\big( \bm X^{\ast}(t)\big)}{t}
 +o(1)\frac{1}{\log n}\int_{1}^n \frac{\big|f\big( \bm X^{\ast}(t)\big)\big|}{t}+o(1) \; a.s.\\
 =& \ep\left[ f\big(\bm X(0)\big)\right]+o(1) \ep\left| f\big(\bm X(0)\big)\right|+o(1)
 \to  \ep\left[ f\big(\bm X(0)\big)\right]=\ep f(\widetilde{\bm C})\;\; a.s. 
 \end{align*}
 \eqref{eq:ASCLT1.1} is proven. Moreover, as shown in Step 2, the exceptional set of probability $0$ can be chosen universally for all such $f$.

{\em Step 4}. We remove the condition $\ep|Z_1|^{2+\epsilon}<\infty$ when   $|f(\bm x)|\le C_0e^{\gamma\|\bm x\bm M\|^2}$.

For this purpose,  we write 
$$ H_n(f;\{\bm x_k\})= \frac{1}{\log n}\sum_{k=1}^n \frac{f\big(\bm x_k/k^{1/2}\big)}{k}  $$
for a sequence $\{\bm x_k\}$ and a function $f$. 
As shown in Step 3,  it is sufficient to   show that 
\begin{enumerate}
  \item[\rm (a)] For $f(\bm x)=(e^{\gamma\|\bm x\bm M\|^2}-A)^+$, 
  \begin{align}\label{eq:expboundness}
 \limsup_{n\to \infty} H_n\left(\big(e^{\gamma \|\bm x\bm M\|^2}-A\big)^+ ;\{(\bm T_k, \bm C_k)\}\right) 
\le    \ep\left[\big(e^{\gamma \|\widetilde{\bm C}\bm M\|^2}-A\big)^+\right];
\end{align}
  \item[\rm (b)] For any bounded Lipschitz function function $f$,
\begin{equation}\label{eq:Lipschitz} \lim_{n\to \infty} H_n\big(f;\{(\bm T_k,\bm C_k)\}\big)=\ep\left[f(\widetilde{\bm C})\right]\;  a.s.   
\end{equation}
\end{enumerate}

We first show that  
\begin{align}
& \label{eq:proofASCLT.3}\limsup_{n\to \infty} H_n((e^{\gamma\|\bm x\|^2}-A)^+;\{(\bm T_k, \bm C_k)\})  \\
\le & \ep\left[\Big(\exp\{8d^2\gamma \sigma_Z^2 (N(0,1))^2\}-A\big)^+\right]+ \ep\left[\Big(\exp\{4\gamma  \mu_Z^2\|   \widetilde{\bm D}\|^2 \}-A\big)^+\right], \nonumber  \\
& \label{eq:proofASCLT.6}  
  \limsup_{n\to \infty} H_n(\|\bm x\|;\{(\bm T_k, \bm C_k)\})  
\le    4d   \sigma_Z  \ep[|N(0,1)|] +2   |\mu_Z| \ep\left[ \|   \widetilde{\bm D}\| \right],
\end{align}
with
$$
 \widetilde{\bm D}\sim N(\bm 0,\bm \Gamma),\; \bm\Gamma=  \begin{pmatrix}    \frac{ 1 }{  1-2\rho } & 
  \frac{1}{ (1-2\rho)(2-\rho)} \\
\frac{1}{ (1-2\rho)(2-\rho)}  &  \frac{2 }{3(1-2\rho)(2-\rho)}\end{pmatrix} \otimes\frac{I_d}{d}.\nonumber
$$
For a sequence $\bm x_k\in\mathbb R^d$, write $g(\bm x_k)=(\bm x_k,\sum_{i=1}^k \bm x_i/k)$. Notice that $f(\bm x)$ is a convex function. Then
$$f(\bm T_k/\sqrt{k}, \bm C_k/\sqrt{k})=f(g(\bm T_k)/\sqrt{k})\le \Big(f(2g(\bm M_{k,2})/\sqrt{k})+f(2\mu_Zg(\bm S_k)/\sqrt{k})\Big)/2, $$
\begin{align*} f(2g(\bm M_{k,2})/\sqrt{k})\le &\big(\exp\{4\gamma\max_{m\le k}\|\bm M_{k,2}\|^2/k\}-A)^+\\
\le & \frac{1}{d}\sum_{j=1}^d \big( \exp\{4\gamma d^2\max_{m\le k}| M_{m,2,j}|^2/k\}-A\big)^+.
\end{align*}

Notice that \eqref{eq:ASCLT1.1} holds for $Z_i\equiv \mu_Z$. Thus, 
$$ \limsup_{n\to \infty}\frac{1}{\log n}\sum_{k=1}^n \frac{f(2\mu_Zg(\bm S_k)/\sqrt{k})}{k}
\le \ep\left[\Big(\exp\{4\gamma  \mu_Z^2\|   \widetilde{\bm D}\|^2 \}-A\Big)^+\right]. $$
When the expectation in the above inequality is finite, the equality holds by \eqref{eq:ASCLT1.1}.  When the expectation is infinite, the inequality is obvious.  

For \eqref{eq:proofASCLT.3}, it is sufficient to show that 
\begin{align}\label{eq:proofASCLT.4} 
 \limsup_{n\to \infty}\frac{1}{\log n}& \sum_{k=1}^n \frac{ \big(\exp\{4\gamma d^2\max_{m\le k}| M_{m,2,j}|^2/k\}-A\big)^+}{k}
\nonumber\\
\le &2\ep\left[\big(\exp\{8d^2\gamma \sigma_Z^2 (N(0,1))^2\}-A\big)^+\right]. 
\end{align}
 
Notice that
$$ \{M_{n,2,j}=\sum_{m=1}^n \sigma_{m,j}(Z_m-\mu_Z); n\ge 1\}\overset{\mathscr{D}}=\{\sum_{m=1}^{\tau_{n,j}}(Z_m-\mu_Z);n\ge 1\} $$
where $\tau_{n,j}=\sum_{m=1}^n \sigma_{m,j}\le n$ (c.f., Doob \cite{Doob1936}, Lemma A.4 of Hu and Zhang \cite{HZ2004}). By the Skorohod embedding theorem, there exist  a standard Brownian motion 
$\{B(t);t\ge 1\}$ and a sequence of i.i.d. non-negative random variables $\{\tau_n;n\ge 1\}$ such that $\ep[\tau_i]=\Var(Z_i)=\sigma_Z^2$ and 
$$ \{\sum_{i=1}^n (Z_i-\mu_Z); n\ge 1\}\overset{\mathscr{D}}=\{B(\sum_{i=1}^n\tau_i);n\ge 1\}. $$ 
On the other hand, 
$$\frac{\sum_{i=1}^n \tau_i}{n}\to \ep[\tau_1]=\sigma_Z^2. $$
Thus, without loss of generality, we can assume that
$$ \max_{m\le k}|M_{m,2,j}|\le \max_{s\le  k }|B(2s\sigma_Z^2)|. $$
For \eqref{eq:proofASCLT.4}, it is sufficient to show that
\begin{align}\label{eq:proofASCLT.5}
  \limsup_{n\to \infty}\frac{1}{\log n}& \sum_{k=1}^n 
\frac{\big( \exp\{\alpha \max_{s\le  k }|B(2s\sigma_Z^2)|^2/k\}-A\big)^+}{k}\nonumber\\
\le & \ep\left[\big(\exp\{\alpha \max_{s\le  1 }|B(2s\sigma_Z^2)|^2\}-A\big)^+\right]. 
\end{align}
Let $X(t)=\max_{s\le  t }|B(2s\sigma_Z^2)|^2/t$, $U_1(t)=\frac{B(2\sigma_Z^2e^t)}{e^{t/2}}$. Then $U_1(t)$ is a Ornstein-Uhlenbeck process which is stationary and ergodic. Then $U(t)=X(e^t)=\sup\limits_{-\infty<y\le 0}\left|U_1(t+y)e^{y/2}\right|$ is a stationary process.  By the Birkhoff   ergodic theorem \cite{Birkhoff1931}, the limit 
\begin{align*}
 &\lim_{n\to \infty}\frac{1}{\log n}\sum_{k=1}^n 
\frac{ \big(\exp\{\alpha \max_{t\le  k }|B(2t\sigma_Z^2)|^2/k\}-A\big)^+}{k}\\
=&\lim_{n\to \infty}\frac{1}{\log n}\int_1^n 
\frac{\big( \exp\{\alpha \max_{s\le  k }|B(2s\sigma_Z^2)|^2/t\}-A\big)^+}{t}dt\\
=&\lim_{n\to \infty}\frac{1}{\log n}\int_0^{\log n} \big(\exp\{\alpha U^2(t)\}-A\big)^+=\overline{f} 
\end{align*}
almost surely exists if $\ep\big[\exp\{\alpha \max\limits_{s\le  1 }|B(2s\sigma_Z^2)|^2\}\big]<\infty$, and also, $$\ep[\overline{f}]=\ep\big[\big(\exp\{\alpha U^2(t)\}-A\big)^+\big]=
\ep\big[\big(\exp\{\alpha \max\limits_{s\le  1 }|B(2s\sigma_Z^2)|^2\}-A\big)^+\big]. $$ 
On the other hand, it is obvious that
$$   \frac{1}{\log n}\sum_{k=1}^n 
\frac{\big( \exp\{\alpha \max_{t_0\le s\le  k }|B(2s\sigma_Z^2)-B(2t_0\sigma_Z^2)|^2/k\}-A\big)^+}{k} $$
has the same limit $\overline{f}$ for any $t_0\ge 0$. Thus, $\overline{f}$ is independent of $B(2t_0\sigma_Z^2)$. Then 
$\overline{f}$ is independent of $\sigma(B(t); t\ge 0)$. Hence, $\overline{f}$ is independent of itself and must be a constant. 
\eqref{eq:proofASCLT.5} is proven, and the proof of \eqref{eq:proofASCLT.3} is completed. With $f(\bm x)=\|\bm x\|$ taking the place of $(e^{\gamma \|\bm x\|^2}-A)^+$, we also have \eqref{eq:proofASCLT.6}. 
 
Now, we begin the proof of \eqref{eq:expboundness}. 
Notice that for any $0<\delta<1$, 
\begin{align*}
 &\big(e^{\gamma \|(\bm x+\bm y)\bm M\|^2}-A\big)^+\\
 \le & \frac{\big(\exp\{\gamma (1+\delta)^2\|\bm x\bm M\|^2\} -A\big)^+}{1+\delta}+ \frac{\delta \big(\exp\{\gamma(1+\delta)^2\|\bm y\bm M\|^2/\delta^2\}-A\big)^+}{1+\delta}\\
 \le & \frac{\big(\exp\{\gamma (1+\delta)^2\|\bm x\bm M\|^2\} -A\big)^+}{1+\delta}+ \frac{\delta\big( \exp\{\alpha\|\bm y\|^2\}-A)^+}{1+\delta}, 
 \end{align*}
 where $\alpha=\gamma(1+\delta)^2c_M/\delta^2$, the fist inequality is due to the convexity of $(e^{\gamma x^2}-A)^+$. 
 
Choose $\delta>0$ small enough such that $\gamma(1+\delta)^3<\gamma_0$.  Denote $Z_i^{(1)}=Z_iI\{|Z_i|\le c\}$, $Z_i^{(2)}=Z_iI\{|Z_i|> c\}$,
$$ \bm T_n^{(j)}=\sum_{k=1}^n \bm \sigma_kZ_k^{(j)}, \; \bm C_n^{(j)}=\frac{1}{n}\sum_{k=1}^n  \bm T_k^{(j)}, \;\; j=1,2. $$
For $(\bm T_n^{(2)}, \bm C_n^{(2)})$, by \eqref{eq:proofASCLT.3} we have
\begin{align*}
&\limsup_{n\to \infty} \frac{1}{\log n}\sum_{k=1}^n\frac{\big(\exp\{\alpha\|(\bm T_k^{(2)}, \bm C_k^{(2)})/\sqrt{k}\|^2/\delta^2\}-A\big)^+}{k} \\
\le &  \ep\left[\big(\exp\big\{8d^2 \alpha  \Var\{Z_1^{(2)}\}  (N(0,1))^2\big\}-A\big)^+\right]
+ \ep\left[\big(\exp\{4 \alpha    (\ep Z_1^{(2)})^2\|   \widetilde{\bm D}\|^2 \}-A\big)^+\right]\\
&\to 2(1-A)^+ =0 \text{ as } c\to \infty, \; \text{ when }\;A>1,
\end{align*}
since $\ep Z_1^{(2)}=\ep[Z_1I\{|Z_1|>c\}]\to 0$ and $\Var\{Z_1^{(2)}\}=\Var(Z_1I\{|Z_1|>c\})\to 0$ as $c\to \infty$. 

Since $Z_i^{(1)}$s are bounded and so satisfies $\ep[ |Z_1^{(1)}|^{2+\epsilon}]<\infty$, for $(\bm T_n^{(1)}, \bm C_n^{(1)})$ and $f(\bm x)=\frac{\big(\exp\{\gamma (1+\delta)^2\|\bm x\bm M\|^2\} -A\big)^+}{1+\delta}$, we have
$$ \limsup_{n\to \infty} \frac{1}{\log n}\sum_{k=1}^n\frac{f\big((\bm T_k^{(1)}, \bm C_k^{(1)})/\sqrt{k}\big)}{k}
\le \ep\left[f(  \widetilde{\bm C}_c)\right] 
$$
as we have shown, where 
$$ \widetilde{\bm C}_c \sim N(\bm 0, \bm\Lambda_c), $$
and $ \bm \Lambda_c$ is defined similar to $\bm \Lambda$ in \eqref{eq:multi-thCMERWCLT1} with $\Var(Z_1I\{|Z_1|\le c\})$, 
$\ep[Z_1I\{|Z_1|\le c\}]$ taking the place of $\sigma_Z^2$, $\mu_Z$, respectively. 
Notice that $ \bm \Lambda_c\to  \bm \Lambda$ as $c\to \infty$. We have
$$ \ep\left[f(  \widetilde{\bm C}_c)\right]\overset{c\to\infty}\longrightarrow \ep\left[f(  \widetilde{\bm C})\right]
\overset{\delta\to 0}\longrightarrow \ep\left[\big(e^{\gamma \|\widetilde{\bm C}\bm M\|^2}-A\big)^+\right]. $$
It follows that \eqref{eq:expboundness} holds. 

For \eqref{eq:Lipschitz},  suppose $|f(\bm x)-f(\bm y)|\le \|\bm x-\bm y\|$. Then
\begin{align*}
&\limsup_{n\to \infty}\left|H_n\big(f;\{(\bm T_k,\bm C_k)\}\big)-H_n\big(f;\{(\bm T_k^{(1)},\bm C_k^{(1)})\}\big)\right|  
\le   \limsup_{n\to \infty}H_n\big(\|\bm x\|;\{(\bm T_k^{(2)},\bm C_k^{(2)})\}\big)\\
\le & 4d   \sqrt{\Var(Z_1I\{|Z_1|> c\})}  \ep[|N(0,1)|] +2   \left|\ep[Z_1I\{|Z_1|> c\}]\right| \ep\left[ \|   \widetilde{\bm D}\| \right] \; a.s.\\
& \to 0 \text{ as } c\to \infty, 
\end{align*}
by \eqref{eq:proofASCLT.6}. And
$$\lim_{n\to \infty} H_n\big(f;\{(\bm T_k^{(1)},\bm C_k^{(1)})\}\big)=\ep\left[f(\widetilde{\bm C}_c)\right]\to \ep\left[f(\widetilde{\bm C})\right] a.s.  \text{ as } c\to \infty. $$
Hence, \eqref{eq:Lipschitz} holds if $f$  is a  bounded Lipschitz function. 
The proof of (i)  for the case of $\rho<1/2$ is now completed. 

Now, suppose $1/2<\rho<1$ and \eqref{eq:non-hom3}. Notice \eqref{eq:processidentical}, \eqref{eq:multi-thRPWInva14} and \eqref{eq:multi-thCMERWInv3}. All the proofs are the same as that for the case of $\rho<1/2$. 

\bigskip
For (ii), we denote $\bm X(t)=\frac{\mu_Z}{\sqrt{d}}(1,2/3)\otimes \widehat{\bm G}_t/{\sqrt{t\log t}}$. Then
$$ \left\{\bm X(t);t> 1\right\}\overset{\mathscr{D}}=\left\{\frac{\mu_Z}{\sqrt{d}}(1,2/3)\otimes \frac{ \bm B(\log t)}{\sqrt{\log t}}; t>1\right\}, $$ 
$$ \left\{\bm U(t)=:\bm X(\exp(e^t));t\in \mathbb R\right\}\overset{\mathscr{D}}=\left\{\frac{\mu_Z}{\sqrt{d}}(1,2/3)\otimes \frac{ \bm B(e^t)}{e^{t/2}}; t\in \mathbb R\right\}. $$ 
Thus, $\bm U(t)$ is an Ornstein-Uhlenbeck process which is stationary and ergodic. By the Birkhoff ergodic theorem, 
\begin{align*}
\lim_{T\to \infty} & \frac{1}{\log\log T}\int_{e}^{T}\frac{f\big(\bm X(t)\big)}{t\log t} dt=  
\lim_{T\to \infty} \frac{1}{\log\log T}\int_{0}^{\log\log T} f\big(\bm U(s)\big)   ds\\
=& \ep\left[f(\bm U(0))\right]=\ep\left[f\big(\frac{\mu_Z}{\sqrt{d}}(1,2/3)\otimes \bm B(1)\big)\right] \; a.s.
\end{align*}
if the expectation is finite. 

By  \eqref{eq:multi-thCMERWInv2}, under $\ep[Z_1^2]<\infty$ we have
$$ \bm T_n =O\big((n\log\log n)^{1/2}\big) + \mu_Z\bm S_n=
\frac{\mu_Z}{\sqrt{d}}\widehat{\bm G}_n+O\big((n\log\log n)^{1/2}\big)\; a.s. $$
and
$$ \bm C_n=\frac{2\mu_Z}{3\sqrt{d}} \widehat{\bm G}_n+O(\sqrt{n\log\log n}) \; a.s.  $$
Thus
$$ \max_{k-1\le t\le k}\left\|\bm X(t)- (\bm T_k,\bm C_k)/\sqrt{k\log k}\right\|= o\big((\log k)^{-1/2+\tau}\big)\; a.s., 0<\tau<1/2. $$
The rest proof is similar to that of \eqref{eq:ASCLTG2} by noticing that $\bm X(t)=O\big((\log\log\log t)^{1/2}\big)$ a.s. 
$\Box$


\bigskip 
  {\bf Acknowledgements.}
 {The author wishes to thank Professor Mikhail Lifshits for enlightening discussions and constructive suggestions, which gave an idea to obtain the precise small ball probabilities of the Gaussian processes considered in this paper.   This work was supported by grants from the NSF of China  (Grant Nos. U23A2064 and 12031005).}


\end{document}